\documentclass[12pt,a4paper]{amsart}
\usepackage{amssymb,amscd}

\theoremstyle{plain}
\newtheorem{theorem}{Theorem}[subsection]
\newtheorem{lemma}[theorem]{Lemma}
\newtheorem{proposition}[theorem]{Proposition}
\newtheorem{corollary}[theorem]{Corollary}

\theoremstyle{definition}

\newtheorem{remark}[theorem]{Remark}
\newtheorem{remarks}[theorem]{Remarks}
\newtheorem{example}[theorem]{Example}

\numberwithin{equation}{subsection}

\newcommand\bC{{\mathbb C}}
\newcommand\bG{{\mathbb G}}
\newcommand\bP{{\mathbb P}}

\newcommand\cE{{\mathcal E}}
\newcommand\cF{{\mathcal F}}
\newcommand\cL{{\mathcal L}}
\newcommand\cM{{\mathcal M}}

\newcommand\cO{{\mathcal O}}
\newcommand\cR{{\mathcal R}}
\newcommand\cS{{\mathcal S}}
\newcommand\cT{{\mathcal T}}
\newcommand\cX{{\mathcal X}}

\newcommand\fa{{\mathfrak a}}

\newcommand\fg{{\mathfrak g}}
\newcommand\fh{{\mathfrak h}}

\newcommand\fU{{\mathfrak U}}

\newcommand\ad{\operatorname{ad}}
\newcommand\aff{\operatorname{aff}}
\newcommand\diag{\operatorname{diag}}
\newcommand\divi{\operatorname{div}}
\newcommand\op{\operatorname{op}}

\newcommand\rk{\operatorname{rk}}

\newcommand\Aut{\operatorname{Aut}}
\newcommand\Ima{\operatorname{Im}}
\newcommand\Ker{\operatorname{Ker}}
\newcommand\Nef{\operatorname{Nef}}
\newcommand\NS{\operatorname{NS}}
\newcommand\Pic{\operatorname{Pic}}
\newcommand\Spec{\operatorname{Spec}}
\newcommand\Supp{\operatorname{Supp}}
\newcommand\Tor{\operatorname{Tor}}

\newcommand\lto{\longrightarrow}

\title[Dolbeault cohomology of log homogeneous varieties]
{Vanishing theorems for Dolbeault \\
cohomology of log homogeneous varieties}

\author{Michel Brion}
\address{Universit\'e de Grenoble I\\
D\'epartement de Math\'ematiques\\
Institut Fourier, UMR 5582 du CNRS\\
38402 Saint-Martin d'H\`eres Cedex, France}
\email{Michel.Brion@ujf-grenoble.fr}

\begin{document}

\begin{abstract}
We consider a complete nonsingular variety $X$ over $\bC$, having a
normal crossing divisor $D$ such that the associated logarithmic
tangent bundle is generated by its global sections. We show that 
$H^i\big(X, L^{-1} \otimes  \Omega_X^j(\log D)\big) = 0$ for any nef
line bundle $L$ on $X$ and all $i < j - c$, where $c$ is an explicit
function of $(X,D,L)$. This implies e.g. the vanishing of 
$H^i(X, L \otimes \Omega_X^j)$ for $L$ ample and $i > j$, and gives
back a vanishing theorem of Broer when $X$ is a flag variety.
\end{abstract}

\maketitle

\section*{Introduction}
\label{sec:introduction}

The main motivation for this work comes from the well-developed
theory of complete intersections in algebraic tori $(\bC^*)^n$ and in
their equivariant compactifications, toric varieties. In particular,
the Hodge numbers of these complete intersections were determined
by Danilov and Khovanskii, and their Hodge structure, by Batyrev, Cox 
and others (see \cite{DK86,BC94,Ma99}). 
This is made possible by the special features of toric geometry; two
key ingredients are the triviality of the logarithmic tangent bundle 
$T_X( - \log D)$, where $X$ is a complete nonsingular toric variety
with boundary $D$, and the Bott--Danilov--Steenbrink vanishing theorem
for Dolbeault cohomology: $H^i(X, L \otimes \Omega_X^j) = 0$ for any
ample line bundle $L$ on $X$ and any $i \geq 1$, $j \geq 0$.

\medskip

A natural problem is to generalise this theory to complete
intersections in algebraic homogeneous spaces and their equivariant
compactifications. As a first observation, the preceding two results
also hold for abelian varieties and, more generally, for the
``semi-abelic'' varieties of \cite{Al02}, that is, equivariant
compactifications of semi-abelian varieties. In fact, for a complete
nonsingular variety $X$ and a divisor $D$ with normal crossings on
$X$, the triviality of the logarithmic tangent bundle is equivalent to
$X$ being semi-abelic with boundary $D$, by a result of Winkelmann
(see \cite{Wi04}). Moreover, it is easy to see that semi-abelic
varieties satisfy Bott vanishing. 

\medskip

The next case to consider after these ``log parallelisable varieties''
should be that of flag varieties. Here counter-examples to Bott
vanishing exist for grassmannians and quadrics, as shown by work of
Snow (see \cite{Sn86}). For example, any  smooth quadric hypersurface  
$X$ in $\bP^{2m}$ satisfies   
$H^{m-1}\big(X, \Omega_X^m(1)\big) \neq 0$. 

\medskip

On the other hand, a vanishing theorem due to Broer asserts that 
$H^i(X, L \otimes \Omega_X^j) = 0$ for any nef line bundle $L$ on a
flag variety $X$, and all $i > j$ (see \cite{Bro93}, and \cite{Bro97}
for a generalisation to homogeneous vector bundles); in this setting,
a line bundle is nef (numerically effective) if and only if it is
effective, or generated by its global sections. Moreover, the same
vanishing theorem holds for any nef line bundle on a complete
simplicial toric variety, in view of a recent result of Mavlyutov 
(see \cite[Thm.~2.4]{Ma08}). 

\medskip

In this article, we obtain generalisations of Broer's vanishing
theorem to any ``log homogeneous'' variety, that is, to a complete
nonsingular variety $X$ having a divisor with normal crossings $D$
such that $T_X(- \log D)$ is generated by its global sections. 
Then $X$ contains only finitely many orbits of
the connected automorphism group $\Aut^0(X,D)$, and these are the strata 
defined by $D$. The class of log homogeneous varieties, introduced and
studied in \cite{Bri07}, contains of course the log parallelisable
varieties, and also the ``wonderful (symmetric) varieties'' of De
Concini--Procesi and Luna (see \cite{DP83, Lu01}). Log homogeneous
varieties are closely related to spherical varieties; in particular,
every spherical homogeneous space has a log homogeneous equivariant
compactification (see \cite{BB96}).

\medskip

Our final result (Theorem \ref{thm:dol}) asserts that 
$H^i(X, L \otimes \Omega_X^j) = 0$ 
\emph{for any nef (resp.~ample) line bundle $L$ on a log
homogeneous variety $X$, and for any $i > j + q + r$ (resp.~$i > j$)}. 
Here $q$ denotes the irregularity of $X$, i.e., the dimension of the
Albanese variety, and $r$ its rank, i.e.,, the codimension of any
closed stratum (these are all isomorphic). Thus, $q + r = 0$ if and
only if $X$ is a flag variety; then our result gives back Broer's
vanishing theorem.

\medskip

We deduce our result from the vanishing of the logarithmic Dolbeault
cohomology groups, 
$H^i\big(X, L^{-1} \otimes \Omega_X^j(\log D)\big)$, for $L$ nef and 
$i < j - c$, where $c \leq q + r$ is an explicit function of $(X,D,L)$; 
see Theorem \ref{thm:vang} for a complete (and optimal) statement. In
particular, $H^i\big(X, \Omega_X^j(\log D)\big) =0$ for all 
$i < j - q - r$; this also holds for all $i > j$, by a general result
on varieties with finitely many orbits (Theorem \ref{thm:finiteg}) .
In view of a logarithmic version of the Lefschetz theorem due to
Norimatsu (see \cite{No78}), this gives information on the mixed Hodge
structure of complete intersections: specifically, 
\emph{for any ample hypersurfaces $Y_1,\ldots,Y_m \subset X$ such that 
$D + Y_1 + \cdots + Y_m$ has normal crossings, the complete
intersection $Y := Y_1 \cap \cdots \cap Y_m$ satisfies 
$H^i\big(Y,\Omega_Y^j(\log D)\big) = 0$ unless 
$i + j \geq \dim(Y)$ or $0 \leq j - i \leq q + r$.}
 
\medskip

Since the proof of our results is somewhat indirect, we first
present it in the setting of flag varieties, and then sketch how to
adapt it to log homogeneous varieties. For a flag variety $X = G/P$,
the tangent bundle $T_X$ is the quotient of the trivial bundle 
$X \times \fg$ (where $\fg$ denotes the Lie algebra of $G$) by the
sub-bundle $R_X$ of isotropy Lie subalgebras. By a homological
argument of ``Koszul duality'' (Lemma \ref{lem:koszul}), Broer's
result is equivalent to the vanishing of $H^i(R_X, p^*L)$ for all 
$i \geq 1$, where $p: R_X \to X$ denotes the structure map. But one
checks that the canonical bundle of the nonsingular variety $R_X$
is trivial, and the projection $f : R_X \to \fg$ is proper, surjective
and generically finite. So the desired vanishing follows from the 
Grauert--Riemenschneider theorem.

\medskip

For an arbitrary log homogeneous variety $X$ with boundary $D$, we
consider the connected algebraic group $G = \Aut^0(X,D)$, with Lie
algebra $\fg = H^0\big(X,T_X(- \log D)\big)$. We may still define
the ``bundle of isotropy Lie subalgebras'' $R_X$ as the kernel of the
(surjective) evaluation map from the trivial bundle $X \times \fg$
to $T_X(-\log D)$, and the resulting map $f : R_X \to \fg$. If $G$ is
linear, we show that the connected components of the 
general fibres of $f$ are toric varieties of dimension $\leq r$
(Theorem \ref{thm:fibres} and Corollary \ref{cor:ineq}, the main
geometric ingredients of the proof). Moreover, any nef line bundle $L$
on $X$ is generated by its global sections. By a generalisation of the 
Grauert--Riemenschneider theorem due to Koll\'ar (see
\cite[Cor.~6.11]{EV92}), it follows that  
$H^i(R_X, p^* L \otimes \omega_{R_X}) = 0$ for any $i > r$. By
homological duality arguments again, this is equivalent to the
vanishing of $H^i\big(X, L^{-1} \otimes \Omega_X^j(\log D)\big)$ for
any such $L$, and all $i < j - r$ (Corollary \ref{cor:van}). In turn,
this easily yields our main result, under the assumption that $G$ is
linear.

\medskip

The case of an arbitrary algebraic group $G$ may be reduced to the
preceding setting, in view of some remarkable properties of the
Albanese morphism of $X$: this is a homogeneous fibration, which
induces a splitting of the logarithmic tangent bundle (Lemma
\ref{lem:split}), and a decomposition of the ample cone (Lemma
\ref{lem:prod}). 

\medskip

The homological arguments of ``Koszul duality'' already appear in the
work of Broer mentioned above, and also in work of Weyman (see
\cite[Chap.~5]{We03}). The latter considers the more general setting
of a sub-bundle of a trivial bundle, but mostly assumes that the
resulting projection is birational, which very seldom holds in our 
setting. 

\medskip

The geometry of the morphism $f: R_X \to \fg$ bears a close analogy
with that of the moment map $\phi : \Omega_X^1 \to \fg^*$, studied
in depth by Knop for a variety $X$ equipped with an action of a
connected reductive group $G$. In particular, Knop considered the
compactified moment map $\Phi : \Omega_X^1( \log D) \to \fg^*$, and he
showed that the connected components of the general fibres of $\Phi$
are toric varieties, if $X$ is log homogeneous under $G$ (see
\cite[p.~265]{Kn94}). By applying another result of Koll\'ar, he also
showed that 
$H^i\big(\Omega_X^1( \log D), \cO_{\Omega_X^1( \log D)}\big) = 0$ 
for all $i\geq 1$ (see \cite[Thm.~4.1]{Kn94}). 

\medskip

However, vanishing results for $H^i\big(\Omega_X^1( \log D), q^*L\big)$,
where $L$ is a nef line bundle on $X$, and $q : \Omega_X^1( \log D) \to X$
denotes the structure map, are only known under retrictive
assumptions on $X$ (see \cite{BB96}). Also, generalising Knop's
vanishing theorem to all log homogeneous varieties (under possibly
non-reductive groups) is an open question.

\medskip

Our construction coincides with that of Knop in the case where $X$ is
a $G \times G$-equivariant compactification of a connected reductive
group $G$: then one may identify $R_X$ with $\Omega_X^1( \log D)$, and
$f$ with the compactified moment map (see Example \ref{exa:gr1}). 
As applications, we obtain very simple descriptions of the algebra of
differential operators on $X$ which preserve $D$, and of the bi-graded
algebra $H^{\bullet}\big(X,\Omega_X^{\bullet}(\log D)\big)$ 
(see Example \ref{exa:gr2}). The structure of the latter algebra also
follows from Deligne's description of the mixed Hodge structure on the
cohomology of $G$ (see \cite[Thm.~9.1.5]{De74}), while the former
seems to be new. In that setting, one may also obtain more precise
information on the numerical invariants of (possibly non-ample)
hypersurfaces in $G$ and $X$; this will be developed elsewhere.

\bigskip

\noindent
{\bf Acknowledgements.}
I wish to thank P.~Brosnan, J-P.~Demailly, L.~Manivel, G.~R\'emond,
and B.~Totaro for stimulating discussions.

\bigskip

\noindent
{\bf Notation.} 
Throughout this article, we consider algebraic schemes, varieties, and
morphisms over the field $\bC$ of complex numbers. We follow the
conventions of the book \cite{Ha77}; in particular, a \emph{variety}
is an integral separated scheme of finite type over $\bC$. By a point,
we always mean a closed point; a \emph{general point} of a variety is
a point of some non-empty Zariski open subset. 

\medskip

We shall consider pairs $(X,D)$, where $X$ is a complete
nonsingular variety of dimension $n$, and $D$ is a 
\emph{simple normal crossing divisor} on $X$, i.e., $D$ is an
effective, reduced divisor with nonsingular irreducible components
$D_1,\ldots,D_{\ell}$ that intersect transversally. We then set
$$
X_0 := X \setminus \Supp(D),
$$ 
the \emph{open part} of $X$. 

\medskip

An \emph{algebraic group} $G$ is a group scheme of finite type over
$\bC$; then each connected component of $G$ is a nonsingular
variety. We denote by $G^0 \subset G$ the \emph{neutral component},
i.e., the connected component through the identity element $e$, and by
$\fg$ the Lie algebra of $G$. 

\medskip

We say that a pair $(X,D)$ is a $G$-\emph{pair}, if $X$ is equipped
with a faithful action of the algebraic group $G$ that preserves $D$.

\section{Logarithmic Dolbeault cohomology of varieties 
with finitely many orbits} 
\label{sec:dc}

\subsection{Differential forms with logarithmic poles}
\label{subsec:dflp}
 
We begin by recalling some basic results on differential forms with
logarithmic poles, referring to \cite[Chap.~2]{EV92} for details. 

Given a pair $(X,D)$ as above and an integer $j \geq 0$, let
$\Omega_X^j(\log D)$ denote the sheaf of \emph{differential forms of
degree $j$ with logarithmic poles along $D$}, consisting of rational
$j$-forms $\omega$ on $X$ such that $\omega$ and $d\omega$ have at
worst simple poles along $D_1,\ldots,D_{\ell}$. The sheaf
$\Omega_X^j(\log D)$ is locally free and satisfies
\begin{equation}\label{eqn:wed}
\Omega_X^j(\log D) = \wedge^j \Omega_X^1(\log D), \quad 
\Omega_X^n(\log D) = \omega_X(D),
\end{equation}
where $\omega_X := \Omega_X^n$ denotes the canonical sheaf. 
As a consequence, the dual sheaf of $\Omega_X^j(\log D)$ is given by 
\begin{equation}\label{eqn:serred}
\Omega_X^j(\log D)^{\vee} = 
\Omega_X^{n - j}(\log D) \otimes \omega_X^{-1}(-D).
\end{equation} 
Moreover,
$$
\Omega_X^1(\log D)^{\vee} =: \cT_X(- \log D)
$$ 
is the \emph{logarithmic tangent sheaf}, i.e., the subsheaf of the
tangent sheaf $\cT_X$ consisting of derivations that preserve the
ideal sheaf of $D$. 

If $(X,D)$ is a $G$-pair for some algebraic group $G$, then the
logarithmic (co)tangent sheaves are $G$-linearised, and we have a 
morphism of linearised sheaves 
\begin{equation}\label{eqn:op}
\op_{X,D} : \cO_X \otimes_{\bC} \fg \lto \cT_X( - \log D),
\end{equation}
the \emph{action map}.

Each divisor $D^k := D - D_k$ induces a simple normal crossing divisor 
on $D_k$, that we denote by $D^k\vert_{D_k}$ or just by $D^k$ for
simplicity. Moreover, taking the residue along $D_k$ yields an exact
sequence
\begin{equation}\label{eqn:resj}
0 \lto \Omega_X^j(\log D^k) \lto \Omega_X^j(\log D) 
\lto \Omega_{D_k}^{j-1}(\log D^k) \lto 0
\end{equation}
that provides an inductive way to relate differential forms with
logarithmic poles to ordinary differential forms. Also, note the exact
sequence 
\begin{equation}\label{eqn:res1}
0 \lto \Omega_X^1 \lto \Omega_X^1(\log D) 
\lto \bigoplus_{k=1}^{\ell} \cO_{D_k} \lto 0.
\end{equation}

Given an invertible sheaf $\cL$ on $X$, we shall consider the 
logarithmic Dolbeault cohomology groups, 
$H^i\big(X, \cL \otimes \Omega_X^j(\log D)\big)$.
Note the isomorphism
\begin{equation}\label{eqn:sd}
H^i\big(X, \cL \otimes \Omega_X^j(\log D)\big)^* \cong 
H^{n-i}\big(X, \cL^{-1}(-D) \otimes \Omega_X^{n-j}(\log D)\big),
\end{equation}
a consequence of Serre duality and (\ref{eqn:serred}). 

Also, recall a vanishing result of Norimatsu (see \cite{No78}): 
\emph{if $\cL$ is ample, then 
$H^i\big(X, \cL^{-1} \otimes \Omega_X^j(\log D)\big) = 0$
for all $i + j < n$.} When $D=0$, this is the classical
Kodaira--Akizuki--Nakano vanishing theorem, and the case of an
arbitrary $D$ follows by induction on the number of irreducible
components of $D$ in view of the exact sequence (\ref{eqn:resj}); see 
\cite[Cor.~6.4]{EV92} for details. 

This vanishing result implies a logarithmic version of the Lefschetz
theorem, also due to Norimatsu (see \cite{No78}). Let $Y_1,\ldots,Y_m$
be ample hypersurfaces in $X$ such that the divisor
$D + Y_1 + \cdots + Y_m$ has simple normal crossings; in particular,
$Y := Y_1 \cap \cdots \cap Y_m$ is a nonsingular complete
intersection, equipped with a simple normal crossing divisor
$D\vert_Y$. Then \emph{the pull-back map
$$
H^i\big(X, \Omega_X^j(\log D)\big) \lto 
H^i\big(Y, \Omega_Y^j(\log D)\big)
$$
is an isomorphism if $i + j < n - m$, and is injective if 
$i + j = n - m$}.

By Hodge theory (see \cite[Sec.~3.2]{De71}), it follows that 
\emph{the pull-back map in cohomology,
$$
H^k(X_0, \bC) \lto H^k(Y_0, \bC),
$$
is an isomorphism for $k < n - m$, and is injective for 
$k = n - m$.}

\subsection{Varieties with finitely many orbits: the linear case}
\label{subsec:vfmo}

In this subsection, we consider a $G$-pair $(X,D)$, where $G$
is a connected algebraic group; we assume that $X$ contains only
finitely many $G$-orbits and that $G$ is linear or, equivalently,
affine. 

We shall obtain a vanishing theorem for the groups 
$H^i\big(X,\Omega_X^j(\log D)\big)$. In the case where $D = 0$, we
have the following result, which is well-known if $G$ is reductive: 

\begin{lemma}\label{lem:pav}
{\rm (i)} With the assumptions of this subsection, $X$ admits a
cellular decomposition (in the sense of \cite[Ex.~1.9.1]{Ful98}). 

\smallskip

\noindent
{\rm (ii)}  There are natural isomorphisms 
$$
A^i(X) \cong H^i(X,\Omega_X^i)
$$ 
for all $i$, where $A^i(X)$ denotes the Chow group of rational
equivalence classes of cycles of codimension $i$, with complex
coefficients. Moreover, 
$$
H^i(X,\Omega_X^j) = 0 \qquad (i \neq j).
$$
\end{lemma}

\begin{proof}
By \cite[Ex.~19.1.11]{Ful98} and Hodge theory, it suffices to show
(i). We shall deduce that assertion from the 
Bia\l ynicki-Birula decomposition (see \cite{Bi73}). 

We may choose a maximal torus $T \subset G$ and a
one-parameter subgroup $\lambda : \bG_m \to T$ such that the fixed
point subscheme $X^{\lambda}$ equals $X^T$. Also, recall that $X^T$ is
nonsingular. For any component $X_i$ of $X^T$, let
$$
X^+_i := \{ x \in X ~\vert~ \lim_{t \to 0} \lambda(t) x \in X_i \}
$$
and 
$$
r_i : X^+_i \lto X_i, \quad 
x \longmapsto \lim_{t \to 0} \lambda(t) x.
$$ 
Then $X$ is the disjoint union of the $X^+_i$, where $X_i$ runs over
the components of $X^T$; moreover, each $X^+_i$ is a locally closed
nonsingular subvariety of $X$, and each retraction $r_i$ is a
locally trivial fibration into affine spaces.

Next, consider the centraliser $G^T \subset G$. Since $G^T$ is
connected and the quotient $G^T/T$ admits no non-trivial subtorus, it
follows that
$$
G^T \cong T \times U,
$$
where $U$ is a unipotent group. Clearly, $U$ acts on $X^T$.

We claim that $X^T$ contains only finitely many orbits of $U$ or,
equivalently, of $G^T$. To check this, it suffices to show that $Z^T$
contains only finitely many $G^T$-orbits, for any $G$-orbit $Z$. 
Given a point $x \in Z^T$, the differential at $e$ of the orbit map
$G \to Z$, $g \mapsto g \cdot x$ yields a surjective map 
$\fg \to T_x Z$ (where $\fg$ denotes the Lie algebra of $G$), and hence
a surjective map $\fg^T \to (T_x Z)^T = T_x(Z^T)$. It follows that
the orbit $G^T \cdot x$ is open in $Z^T$, which implies our claim.

Note that $U$ preserves each component $X_i$, and $r_i$ is
$U$-equivariant. Thus, given an orbit $Z = U \cdot x$ in $X_i$, the
preimage $r_i^{-1}(Z)$ is isomorphic to the homogeneous bundle
$U \times^{U_x} F$ over $Z \cong U/U_x$, where $U_x$ denotes the
isotropy group of $x$ in $U$, and $F$ the fibre of $r_i$ at $x$ 
(so that $F$ is preserved by $U_x$ and isomorphic to an affine
space). Since $U$ is unipotent, it contains a closed subvariety $V$,
isomorphic to an affine space, such that the multiplication of the
group $U$ induces an isomorphism $V \times U_x \cong U$. Then 
$r_i^{-1}(Z) \cong V \times F$ is an affine space, which yields the
desired cellular decomposition.
\end{proof}

From that result and inductive arguments, we shall deduce:

\begin{theorem}\label{thm:finite}
With the assumptions of this subsection, there are natural
isomorphisms
\begin{equation}\label{eqn:chow}
A^i(X_0) \cong H^i\big(X, \Omega_X^i(\log D)\big)
\end{equation}
for all $i$. Moreover, 
\begin{equation}\label{eqn:diag}
H^i\big(X, \Omega_X^j(\log D)\big) = 0 \qquad (i > j).
\end{equation}
\end{theorem}

\begin{proof}
We may assume that $D \neq 0$ in view of Lemma \ref{lem:pav}.
We first prove the vanishing assertion (\ref{eqn:diag}), by induction
on the number of irreducible components of $D$ and the dimension of
$X$. Write $D = D_1 + D^1$, where $D_1$ is irreducible. Then
(\ref{eqn:resj}) yields an exact sequence
$$
H^i\big(X, \Omega_X^j(\log D^1)\big) \lto
H^i\big(X, \Omega_X^j(\log D)\big) \lto
H^i\big(D_1, \Omega_{D_1}^{j-1}(\log D^1)\big)
$$
which implies our assertion.

Next, we construct the isomorphism (\ref{eqn:chow}).
If $D$ is irreducible, then the exact sequence
$0 \to \Omega_X^i \to \Omega_X^i(\log D) \to \Omega_D^{i-1} \to 0$
(a special case of (\ref{eqn:resj})) yields a diagram
$$
\CD
A^{i-1}(D) @>>> A^i(X) @>>> A^i(X \setminus D) @>>> 0 \\
@VVV @VVV \\
H^{i-1}(D, \Omega_D^{i-1}) @>>> H^i(X, \Omega_X^i) @>>> 
H^i\big(X, \Omega_X^i(\log D)\big) @>>> 0, \\
\endCD
$$
where the top row is a standard exact sequence of Chow groups
(see \cite[Prop.~1.8]{Ful98}), the bottom row is exact since
$H^i(D, \Omega_D^{i-1})=0$, the vertical arrows are isomorphisms, and
the square commutes by functoriality of the cycle map (see
\cite[Sec.~19.1]{Ful98}). This yields the desired isomorphism.

In the general case, we argue again by induction on the number of
irreducible components of $D$ and the dimension of $X$. The exact
sequences
$$
H^{i-1}(D_k, \Omega_{D_k}^{i-1}) \lto H^i(X, \Omega_X^i) \lto 
H^i\big(X, \Omega_X^i(\log D_k)\big) \lto 0
$$
for $k = 1, \ldots, \ell$ and the natural maps
$$ 
H^i\big(X, \Omega_X^i(\log D_k)\big) \lto
H^i\big(X, \Omega_X^i(\log D)\big)
$$
yield a complex
\begin{equation}\label{eqn:comp}
\bigoplus_{k=1}^{\ell} H^{i-1}(D_k, \Omega_{D_k}^{i-1}) \lto 
H^i(X, \Omega_X^i) \lto H^i\big(X, \Omega_X^i(\log D)\big) \lto 0.
\end{equation}
We claim that this complex is exact. Consider indeed the
commutative diagram
$$
\CD
H^{i-1}(D_1, \Omega_{D_1}^{i-1}) &
\to & H^i(X, \Omega_X^i) &
\to & H^i\big(X, \Omega_X^i(\log D_1)\big)  & \to 0 \\
@VVV   @VVV    @VVV  \\
H^{i-1}\big(D_1, \Omega_{D_1}^{i-1}(\log D^1)\big) &
\to & H^i\big(X, \Omega_X^i (\log D^1)\big) &
\to & H^i\big(X, \Omega_X^i(\log D)\big) & \to 0. \\
\endCD
$$
Then the rows are exact, by the preceding argument and the vanishing of 
$H^i\big(D_1, \Omega_{D_1}^{i-1}(\log D^1)\big)$. Moreover, the left and
middle vertical maps are surjective by the induction assumption; thus,
so is the right vertical map. This implies our claim.

That claim implies in turn the isomorphism (\ref{eqn:chow}), by
comparing (\ref{eqn:comp}) with the complex
$$
\bigoplus_{k=1}^{\ell} A^{i-1}(D_k) \lto A^i(X) 
\lto A^i(X_0) \lto 0
$$
which is exact in view of the standard exact sequence
$$
A^{i-1}\big(\Supp(D)\big) \lto 
A^i(X) \lto A^i(X_0) \lto 0
$$
and the surjectivity of the natural map 
$$
\bigoplus_{k=1}^{\ell} A^{i-1}(D_k) \lto
A^{i-1}\big(\Supp(D)\big).
$$
\end{proof}

\begin{remark}
One can show that $X_0$ is a linear variety, as defined by Totaro 
in \cite{To08}. By Theorem 3 of that article, it follows that the
Chow group of $X_0$ is isomorphic to the smallest subspace of 
Borel-Moore homology (with complex coefficients) with respect to the
weight filtration. In turn, this yields another proof of Theorem
\ref{thm:finite}, admittedly less direct than the proof presented here.
\end{remark}

\subsection{Varieties with finitely many orbits: the general case}
\label{subsec:vfmog}

We still consider a $G$-pair $(X,D)$, where $G$ is a connected
algebraic group, and $X$ contains only finitely many $G$-orbits, but
we no longer assume that $G$ is linear.

We shall obtain a generalisation of Theorem \ref{thm:finite} to that
setting; for this, we recall (after \cite[Prop.~2.4.1]{Bri07}) a
reduction to the linear case, via the Albanese morphism
$$
\alpha: X \lto A
$$
(the universal morphism to an abelian variety).

By Chevalley's structure theorem, $G$ admits a largest connected
affine subgroup $G_{\aff}$; this subgroup is normal in $G$, and the
quotient $G/G_{\aff}$ is an abelian variety. Moreover, $X$ is
equivariantly isomorphic to the total space of a homogeneous bundle 
$G \times^I Y$, where $I\subset G$ is a closed subgroup with neutral
component $G_{\aff}$, and $Y \subset X$ is a closed nonsingular
subvariety, preserved by $I$; both $I$ and $Y$ are unique. 

As a consequence, $I$ is a normal affine subgroup of $G$, and the
quotient $G/I$ is an abelian variety equipped with an isogeny
$G/G_{\aff} \to G/I$. Moreover, $I$ acts on $Y$ with
finitely many orbits, and $D$ induces a simple normal crossing divisor 
$E$ on $Y$, preserved by $I$. Finally, the Albanese morphism
$\alpha$ may be identified to the homogeneous fibration 
$G \times^I Y \to G/I$ with fibre $Y$.

Also, note a splitting property of the logarithmic (co)tangent sheaf:

\begin{lemma}\label{lem:split}
With the preceding notations, there is an isomorphism of
$G$-linearised sheaves
\begin{equation}\label{eqn:split}
\Omega_X^1(\log D) \cong 
\Omega_{X/A}^1( \log D) \oplus (\cO_X \otimes_{\bC} \fa^*),
\end{equation}
where $\Omega_{X/A}^1( \log D)\otimes \cO_Y \cong \Omega_Y^1(\log E)$,
and $\fa$ denotes the Lie algebra of $A$ (so that $G$ acts trivially
on $\fa$).
 
Moreover, the composite map 
$\cO_X \otimes_{\bC} \fa^* \to \Omega_X^1(\log D) 
\to \cO_X \otimes_{\bC} \fg^*$
(where the map on the right is the transpose of the action map
(\ref{eqn:op})) is induced from the map $\fa^* \to \fg^*$, the transpose 
of the quotient map $\fg \to \fg/\fg_{\aff} = \fa$.
\end{lemma}

\begin{proof}
By \cite[Prop.~2.4.1]{Bri07}, the Albanese fibration yields an
exact sequence of $G$-linearised sheaves
\begin{equation}\label{eqn:alb}
0 \lto \cT_{X/A}(- \log D) \lto \cT_X(- \log D) \lto \alpha^*\cT_A \lto 0,
\end{equation}
where $\cT_{X/A}(- \log D)\otimes \cO_Y \cong \cT_Y( -\log E)$. 
Also, note that $\alpha^*\cT_A \cong \cO_X \otimes \fa$. Since
$G$ acts on $A$, and $\alpha$ is equivariant, the composite map 
$$
\CD
\cO_X \otimes_{\bC} \fg @>{\op_{X,D}}>> \cT_X(- \log D) 
@>>> \cO_X \otimes_{\bC} \fa
\endCD
$$
is induced from the quotient map $\fg \to \fa$.

Choose a subspace $\tilde{\fa}\subset \fg$ such that the composite map
$\tilde{\fa} \to \fg  \to \fa$ is an isomorphism. Then the composite map
$$
\cO_X \otimes_{\bC} \tilde{\fa} \lto \cO_X \otimes_{\bC} \fg 
\lto \cT_X(- \log D) \lto \cO_X \otimes_{\bC} \fa
$$
is an isomorphism as well; thus, the exact sequence (\ref{eqn:alb}) is
split. Taking duals yields our assertions.
\end{proof}

\begin{remark}
With the notation of the preceding proof, we may further assume that
$\tilde{\fa}$ is contained in the centre of $\fg$. Indeed, if $C(G)$
denotes the centre of the group $G$, then the natural map $C(G) \to A$
is surjective, as follows e.g. from \cite[Lem.~1.1.1]{Bri07}. 

This yields a decomposition of the logarithmic tangent sheaf into a
direct sum of the integrable sub-sheaves $\cT_{X/A}(- \log D)$ and
$\cO_X \otimes_{\bC} \tilde{\fa}$, and an analogous splitting of the
tangent sheaf.
\end{remark}

We now come to the main result of this section:

\begin{theorem}\label{thm:finiteg}
With the notation and assumptions of this subsection, 
\begin{equation}\label{eqn:diagg}
H^i\big(X, \Omega_X^j(\log D)\big) = 0 \qquad \big(i > j + q(X)\big),
\end{equation}
where $q(X):= \dim(A)$ denotes the irregularity of $X$. Moreover,
there is an isomorphism
\begin{equation}\label{eqn:chowg}
H^{j + q(X)}\big(X, \Omega_X^j(\log D)\big) \cong A^j(Y_0)^I
\end{equation}
(the subspace of $I$-invariants in $A^j(Y_0)$), and $I$ acts in
$A^j(Y_0)$ via the finite quotient $I/G_{\aff}$. 
\end{theorem}

\begin{proof}
Lemma \ref{lem:split} yields decompositions
\begin{equation}\label{eqn:albj}
\Omega_X^j(\log D) \cong \bigoplus_{k=0}^j 
\Omega_{X/A}^k(\log D) \otimes_{\bC} \wedge^{j-k} \fa^*
\end{equation}
and isomorphisms
$$
\Omega_{X/A}^k(\log D)\otimes \cO_Y \cong \Omega_Y^k(\log E).
$$
Since $H^i\big(Y, \Omega_Y^k(\log E)\big) = 0$ for $i > k$ by 
(\ref{eqn:diag}), this yields in turn
\begin{equation}
R^i \alpha_* \Omega_{X/A}^k(\log D)= 0 \qquad (i > k).
\end{equation}
Together with (\ref{eqn:albj}), it follows that 
\begin{equation}\label{eqn:hdi}
R^i \alpha_* \Omega_X^j(\log D)= 0 \qquad (i > j).
\end{equation}
This implies the vanishing (\ref{eqn:diagg}) in view of the Leray
spectral sequence
$$
H^p\big(A, R^q \alpha_* \Omega_X^j(\log D)\big) \Rightarrow
H^{p+q}\big(X, \Omega_X^j(\log D)\big)
$$
which also yields isomorphisms
$$
H^{q(X) + j}\big(X, \Omega_X^j(\log D)\big) \cong 
H^{q(X)}\big(A, R^j \alpha_* \Omega_X^j(\log D)\big).
$$
Moreover, 
$R^j \alpha_* \Omega_X^j(\log D ) = 
R^j \alpha_* \Omega_{X/A}^j(\log D )$ 
is the $G$-linearised sheaf on $A = G/I$ associated with the $I$-module
$H^j\big(Y, \Omega_Y^j(\log E)\big)$, i.e., with $A^j(Y_0)$ in
view of (\ref{eqn:chow}). By Serre duality, it follows that
$$
H^{q(X)}\big(A, R^j \alpha_* \Omega_X^j(\log D)\big) \cong H^0(A, \cF)^*,
$$
where $\cF$ denotes the $G$-linearised sheaf on $G/I$ associated with
the dual $I$-module $A^j(Y_0)^*$. Since the connected linear algebraic
group $G_{\aff}$ acts trivially on the Chow group $A^j(Y_0)$, we have
$$\displaylines{
H^0(A, \cF) = H^0(G/I, \cF) \cong 
\big(\cO(G) \otimes_{\bC} A^j(Y_0)^*\big)^I
\hfill \cr \hfill
\cong 
\big(\cO(G)^{G_{\aff}} \otimes_{\bC} A^j(Y_0)^*\big)^{I/G_{\aff}}
\cong \big(A^j(Y_0)^*\big)^{I/G_{\aff}}.
\cr}$$
This yields the isomorphism (\ref{eqn:chowg}).
\end{proof}

\begin{remarks}\label{rem:cs}
(i) Let $Y_1,\ldots, Y_m \subset X$ be ample hypersurfaces such that
the divisor $D + Y_1 + \cdots + Y_m$ has simple normal crossings, and
consider the complete intersection $Y := Y_1 \cap \cdots \cap Y_m$. 
Combining Theorem \ref{thm:finiteg} with the logarithmic Lefschetz
theorem recalled in Subsec.~\ref{subsec:dflp}, we see that
$$
H^i\big(Y, \Omega_Y^j(\log D) \big) = 0 
\qquad \big (i + j < n - m \text{ and } i > j + q(X)\big).
$$

If $G$ is linear, i.e., $q(X) = 0$, then also
$$
H^i\big(Y, \Omega_Y^i(\log D) \big) \cong A^i(X_0) \qquad
(2i < n - m).
$$ 
Moreover, if $i = \frac{n - m}{2}$, then
$A^i(X_0) \hookrightarrow H^i\big(Y, \Omega_Y^i(\log D) \big)$. 
In particular, 
$A^i(X_0) \hookrightarrow A^i(Y_0)$ for all $i \leq \frac{n - m}{2}$.

If, in addition, $X_0$ is affine, then 
$H^i\big(Y, \Omega_Y^j(\log D)\big) = 0$ whenever  
$i + j > n - m$. Thus, the only ``unknown'' groups
$H^i\big(Y, \Omega_Y^j(\log D) \big)$ are those where $i + j = n - m$.
To compute their dimension, it suffices to determine the Euler
characteristic $\chi\big(Y, \Omega_Y^j(\log D) \big)$, which is
expressed in topological terms via the Riemann--Roch theorem. 

Since the topological Euler characteristic of $Y_0$
satisfies an adjunction formula due to Norimatsu and Kiritchenko 
(see \cite{No78, Ki06}), this yields a determination of the Betti 
numbers of $Y_0$, as already observed in \cite{DK86} for toric
varieties. 

\smallskip

\noindent
(ii) Taking $D = 0$ in Theorem \ref{thm:finiteg} and using Serre
duality yields the vanishing
\begin{equation}\label{eqn:dia}
H^i(X, \Omega_X^j) = 0 \qquad \big(\vert i - j \vert > q(X)\big)
\end{equation}
for a complete nonsingular variety $X$ where an algebraic group $G$
acts with finitely many orbits.

This is closely related to a result of Carrell and Lieberman (see
\cite[Thm.~1]{CL73}): 
\begin{equation}\label{eqn:cl}
H^i(M, \Omega_M^j) = 0 \qquad \big(\vert i - j \vert > \dim Z(V)\big),
\end{equation} 
where $M$ is a compact K\"ahler manifold admitting a global vector
field $V$ with non-empty scheme of zeros $Z(V)$. 

In fact, (\ref{eqn:cl}) implies (\ref{eqn:dia}) when $X$ is projective
and the maximal connected affine subgroup $G_{\aff} \subset G$ is
reductive. Consider indeed a general one-parameter subgroup 
$\lambda : \bG_m \to G_{\aff}$ and the associated vector
field $V \in \fg_{\aff}$. Then $Z(V) = X^{\lambda} = X^T$, where $T$
denotes the unique maximal torus of $G_{\aff}$ containing the image of
$\lambda$. It follows that $Z(V)$ meets each $G_{\aff}$-orbit along a
finite set, non-empty if the orbit is closed; as a consequence,
$\dim Z(V) = q(X)$. However, it is not clear whether (\ref{eqn:dia})
may be deduced from (\ref{eqn:cl}) when (say) $G_{\aff}$ is
unipotent.
\end{remarks}

\section{Geometry of log homogeneous varieties} 
\label{sec:lhv}

\subsection{Basic properties}
\label{subsec:lhv}
Let $(X,D)$ be a $G$-pair, where $G$ is a connected algebraic group. 
Following \cite{Bri07}, we say that $X$ is \emph{log homogeneous under
$G$ with boundary $D$}, or that $(X,D)$ is $G$-\emph{homogeneous},
if the action map (\ref{eqn:op}) is surjective. 

Under that assumption, the open part $X_0$ is a unique $G$-orbit. 
Moreover, by \cite[Cor.~3.2.2]{Bri07}, the $G$-orbit closures in $X$
are exactly the non-empty partial intersections of boundary divisors,
$$
D_{k_1, \ldots, k_m} := D_{k_1} \cap \ldots \cap D_{k_m},
$$
and each $D_{k_1, \ldots, k_m}$ is log homogeneous under $G$ with
boundary being the restriction of the divisor
$$
D^{k_1, \ldots, k_m} := 
\sum_{k \neq k_1,\ldots,k_m} D_k.
$$
In particular, $X$ contains only finitely many $G$-orbits, and these
coincide with the orbits of the connected automorphism group
$\Aut^0(X,D)$. Moreover, the closed orbits are exactly the minimal
non-empty partial intersections. By \cite[Thm.~3.3.3]{Bri07}, these
closed orbits are all isomorphic. We call their common codimension the 
\emph{rank} of $X$, and denote it by $\rk(X)$. Note that 
\begin{equation}\label{eqn:rk}
\rk(D_{k_1,\ldots,k_m}) = \rk(X) - m.
\end{equation}

Log homogeneity is preserved by equivariant blowing up, in the
following sense. Let $X'$ be a complete nonsingular $G$-variety
equipped with a $G$-equivariant birational morphism 
$u : X' \to X$. Denote by $D'$ the reduced inverse image of $D$. 
Then $(X',D')$ is a homogeneous $G$-pair, by
\cite[Prop.~2.3.2]{Bri07}. We now show that logarithmic Dolbeault
cohomology is also preserved:

\begin{lemma}\label{lem:ind}
With the preceding notation and assumptions, there are isomorphisms
\begin{equation}\label{eqn:iso}
u^*\Omega_X^j(\log D) \cong \Omega_{X'}^j(\log D').
\end{equation}
Moreover, any invertible sheaf $\cL$ on $X$ satisfies
\begin{equation}\label{eqn:ind}
H^i\big(X', u^*\cL \otimes \Omega_{X'}^j(\log D')\big) \cong
H^i\big(X, \cL \otimes \Omega_X^j(\log D)\big)
\end{equation}
for all $i$ and $j$.
\end{lemma}

\begin{proof}
The natural morphism
$$
du: \cT_{X'}(- \log D') \lto u^* \cT_X(- \log D)
$$
is clearly surjective, and hence is an isomorphism since its source and
target are locally free sheaves of the same rank. This implies
(\ref{eqn:iso}) and, in turn, the isomorphism (\ref{eqn:ind}) by using
the projection formula and the equalities $u_* \cO_{X'} = \cO_X$, 
$R^i u_* \cO_{X'} = 0$ for all $i \geq 1$. 
\end{proof}

\subsection{The bundle of isotropy Lie subalgebras}
\label{subsec:bils}

We still consider a pair $(X,D)$, homogeneous under a connected
algebraic group $G$. The action map (\ref{eqn:op}) yields an exact
sequence
\begin{equation}\label{eqn:exa}
0 \lto \cR_X \lto \cO_X \otimes_{\bC} \fg \lto \cT_X(-\log D) \lto 0,
\end{equation}
where $\cR_X$ is locally free; equivalently, we have an exact sequence
\begin{equation}\label{eqn:dual}
0 \lto \Omega_X^1(\log D) \lto \cO_X \otimes_{\bC} \fg^* \lto \cR_X^{\vee}
\lto 0. 
\end{equation}

We denote by $R_X$ the vector bundle over $X$ associated with the
locally free sheaf $\cR_X$. Specifically, the structure map
$$
p: R_X \lto X
$$ 
satisfies
\begin{equation}\label{eqn:dir}
p_* \cO_{R_X} = \bigoplus_{j \geq 0} S^j \cR_X^{\vee},
\end{equation}
where $S^j$ denotes the $j$-th symmetric power over $\cO_X$.

By \cite[Prop.~2.1.2]{Bri07}, the fibre $R_{X,x}$ at an arbitrary 
point $x$ is an ideal of the isotropy Lie subalgebra $\fg_x$:
the kernel $\fg_{(x)}$ of the representation of $\fg_x$ in the
normal space to the orbit $G \cdot x$ at $x$. In particular, 
$R_{X,x} = \fg_x$ if $x \in X_0$. We may thus call $R_X$ the 
\emph{bundle of isotropy Lie subalgebras}.

We may view $R_X$ as a closed $G$-stable subvariety of $X \times \fg$,
and denote by
$$
f: R_X \lto \fg
$$
the second projection. Then $f$ is proper, $G$-equivariant, and its
fibres may be identified to closed subschemes of $X$ via the first
projection $p$. 

Since $X$ is complete, the vector bundle $R_X$ is trivial if and only
if $f$ is constant, i.e., $\fg_{(x)}$ is independent of $x \in X$.  By 
\cite[Thm.~2.5.1]{Bri07}, this is also equivalent to $X$ being a
semi-abelic variety.

Returning to an arbitrary homogeneous $G$-pair $(X,D)$, choose a point
$x_0 \in X_0$ and denote by 
$$
H := G_{x_0}
$$ 
its isotropy group, with Lie algebra 
$$
\fh := \fg_{x_0}.
$$ 
This identifies $X_0$ to the homogeneous space $G/H$, and the 
pull-back $R_{X_0}$ to the homogeneous vector bundle
$G \times^H \fh$ associated with the adjoint representation of 
the (possibly non-connected) algebraic group $H$. The
restriction 
$$
f_0: R_{X_0} \lto \fg
$$ 
is identified to the ``collapsing'' morphism
\begin{equation}\label{eqn:col}
G \times^H \fh \lto \fg, 
\quad (g,\xi) H \longmapsto g \cdot \xi
\end{equation}
where the dot denotes the adjoint action. In particular,
\begin{equation}\label{eqn:ima}
f(R_X) = \overline{G \cdot \fh}
\end{equation}
and for any $\xi \in \fh$, the fibre of $f_0$ at the point 
$(1,\xi)H \in R_{X,x_0}$ is identified to the fixed point subscheme
$(G/H)^{\xi}$.

\subsection{The general fixed point subschemes}
\label{subsec:gfps}

We keep the notation and assumptions of Subsec.~\ref{subsec:bils}, and
assume in addition that the algebraic group $G$ is linear.

\begin{theorem}\label{thm:fibres}
With the preceding assumptions, the connected components of the
general fibres of $f$ are toric varieties under subtori of $G$, of
dimension $\rk(G) - \rk(H)$.
\end{theorem} 

\begin{proof}
It suffices to consider fibres at points $(1,\xi) H$, where $\xi$ is a 
general point of $\fh$. Let $\xi = s + n$ be the Jordan decomposition,
that is, $s \in \fh$ is semi-simple, $n \in \fh$ is nilpotent, and
$[s,n] = 0$. Then $\rk(G) = \rk(G^s)$, where $G^s$ denotes the
centralizer of $s$ in $G$; likewise, $\rk(H) = \rk(H^s)$. Also, note
that $n \in \fh^s$ (the Lie algebra of $H^s$). 
Moreover, since $\xi$ is general, we may assume that $s$ is a general
point of the Lie algebra of a maximal torus $T_H \subset H$. Then
$G^s = G^{T_H}$ acts on the fixed point subscheme $X^s = X^{T_H}$
through the quotient group $G^{T_H}/T_H$; moreover, $(H^0)^{T_H}/T_H$
is unipotent. Together with Lemma \ref{lem:loc}, this yields a
reduction to the case where $\rk(H) = 0$; equivalently, $H^0$ is
unipotent.

Under that assumption, we claim that $H^0$ is a maximal unipotent
subgroup of $G$. Indeed, the variety $G/H$ is spherical under
any Levi subgroup $L$ of $G$, by \cite[Thm.~3.2.1]{Bri07}. Therefore, 
$\dim(G/H) \leq \dim(B_L)$, where $B_L$ denotes a Borel subgroup of
$L$. But $\dim(B_L) = \dim(G/U)$, where $U \subset G$ is a maximal
unipotent subgroup. Thus, $\dim(H) \geq \dim(U)$ which implies our claim.

By that claim, the normalizer of $H$ is a Borel subgroup of $G$,
that we denote by $B$. Moreover, for any maximal torus $T \subset B$,
the intersection $T \cap H$ is finite, and $B/H \cong T/(T\cap H)$. 
Since $B$ normalizes $\fh$, the morphism 
$f_0 : G \times^H \fh \to \fg$ factors as the natural map
$$
u : G \times^H \fh \lto G \times^B \fh
$$
followed by the collapsing morphism
$$
v : G \times^B \fh \lto \fg.
$$
Also, $v$ is birational onto its image, the cone of nilpotent elements
in $\fg$ (this is well-known in the case where $G$ is reductive,
and the general case follows by using a Levi decomposition of $G$). On
the other hand, the fibre of $u$ at any point $(1,\xi)H$ is $B/H$,
where both are identified to subvarieties of $G/H$. This identifies
the general fibres of $f_0$ to quotients of maximal tori of $G$ by
finite subgroups.
\end{proof}

\begin{lemma}\label{lem:loc}
Let $s$ be a semi-simple element of $\fh$ with centralizer $G^s$ and
fixed point subscheme $X^s$, and let $Y$ denote the connected
component of $X^s$ through $x_0$. Then:

\smallskip

\noindent
{\rm (i)} $Y$ is a log homogeneous variety under $G^s$ with boundary
$D \vert_Y$.  Moreover, the fixed point subscheme $R_X^s$ is a vector
bundle over $X^s$, and its pull-back to $Y$ is $R_Y$.

\smallskip

\noindent
{\rm (ii)} For any nilpotent element $n \in \fh^s$, the fibres of  
$f : R_X \to \fg$ and of $f_Y : R_Y \to \fg^s$ at $s + n \in \fh^s$
coincide in a neighborhood of the point $(1, s + n) H$.
\end{lemma}

\begin{proof}
(i) Since $s$ is semi-simple, $X^s$ is nonsingular and 
$T_x(X^s) = (T_x X) ^s$ for any $x \in X^s$. Moreover, 
$T_x\big((G \cdot x)^s\big) = \big(T_x (G\cdot x)\big)^s 
= \fg^s \cdot x$. 
As a consequence, $G^s \cdot x$ is a component of $(G\cdot x)^s$. It
follows readily that $D$ induces a divisor with normal crossings on
$X^s$. Moreover, the exact sequence
$$ 
0 \lto \fg_{(x)} \lto \fg \lto T_x X (- \log D) \lto 0
$$
(see \cite[(2.1.4)]{Bri07}) yields an exact sequence 
$$ 
0 \lto \fg_{(x)}^s \lto \fg^s \lto T_x X^s (- \log D) \lto 0.
$$
This implies both assertions.

(ii) It suffices to show that both fibres have the same tangent
space at the point $(1,\xi)H$, where $\xi := s + n$. For this,
consider the map
$$
\Phi : G \times \fh \lto \fg, \quad 
(g,\xi) \longmapsto g \cdot \xi,
$$
invariant under the action of $H$ via 
$h \cdot (g,\xi) = (gh^{-1}, h \cdot \xi)$
and which induces the map $f_0$ on the quotient $G \times^H \fh$. The
differential of $\Phi$ at $(1,\xi) H$ may be identified to the map
$$
\varphi : \fg \times \fh \lto \fg, \quad 
(u,v) \longmapsto [u,\xi] + v
$$
and the differential of the orbit map 
$$
H \lto G \times \fh, \quad 
h \longmapsto (h^{-1}, h \cdot \xi)
$$
is identified with the map
$$
\psi: \fh \lto \fg \times \fh, \quad 
w \longmapsto (-w, [w,\xi]).
$$
Thus, the tangent space at $(1,\xi) H$ of the fibre of $f$ through
that point is the homology space $\Ker(\varphi)/\Ima(\psi)$.

Next, consider the decomposition 
$\fg = \bigoplus_{\lambda \in \bC} \fg_{\lambda}$ into eigenspaces of
$s$, where $\fg_0 = \fg^s$, and the induced decomposition 
$\fh = \bigoplus_{\lambda \in \bC} \fh_{\lambda}$. For 
$(u,v) \in \Ker(\varphi)$, this yields with obvious notation:
$$
v_\lambda - \big(\lambda + \ad(n)\big) u_\lambda =0 \quad
\text{ for all }\lambda.
$$
If $\lambda \neq 0$, then $\lambda + \ad(n)$ is an automorphism of
$\fg$ preserving $\fh$, and hence 
$(u_\lambda,v_\lambda) \in \Ima(\psi)$. Thus, 
$$
\Ker(\varphi)/\Ima(\psi) = \Ker(\varphi^s)/\Ima(\psi^s)
$$
with obvious notation again; this yields the desired equality of
tangent spaces.
\end{proof}

\begin{remark}\label{rem:red}
Given an arbitrary closed subgroup $H$ of a connected linear
algebraic group $G$, and a semi-simple element $s \in \fh$, the orbit 
$G^s \cdot x_0 \cong G^s/H^s$ is open in the fixed point subscheme
$(G/H)^s$. If $G$ and $H$ are reductive, then any general point 
$s \in \fh$ is semi-simple, and $(H^s)^0$ is a maximal torus of
$H$, contained in the centre of the connected reductive group
$G^s$. If, in addition, the homogeneous space $G/H$ is spherical, then
$G^s/H^s$ is also spherical, as follows e.g. from Lemma
\ref{lem:loc}. This implies that $G^s$ is a torus, and yields a very
simple proof of Theorem \ref{thm:fibres} in that special case. 
\end{remark}

\begin{example}\label{exa:gr1}
Any connected reductive group $G$ may be viewed as the
homogeneous space $(G \times G)/ \diag(G)$ for the action of 
$G \times G$ by left and right multiplication. This homogeneous space
is spherical, and hence admits a log homogeneous equivariant
compactification $X$. 

We claim that 
\begin{equation}\label{eqn:grp}
\cR_X \cong \Omega_X^1(\log D)
\end{equation}
as $G \times G$-linearised sheaves; in other words, $R_X$ is
equivariantly isomorphic to the total space of the logarithmic
cotangent bundle. To see this, choose a non-degenerate $G$-invariant
quadratic form $q$ on $\fg$; then the quadratic form 
$(q, -q)$ on $\fg \times \fg$ is non-degenerate and 
$G \times G$-invariant. Moreover, the fibre of $R_X$ at the identity
element, $\fg_{(e)} = \diag(\fg)$, is a Lagrangian subspace of 
$\fg \times \fg$. It follows that $R_X$ is a Lagrangian sub-bundle of
the trivial bundle $X \times \fg \times \fg$ equipped with the
quadratic form $(q,-q)$. Thus, the quotient bundle 
$(X \times \fg \times \fg)/R_X$ is isomorphic to the dual of
$R_X$. This implies our claim in view of the exact sequence
(\ref{eqn:exa}).

In fact, via the isomorphism $\fg \times \fg \cong \fg^* \times \fg^*$
defined by $(q,-q)$, the map $f: R_X \to \fg \times \fg$ is identified to
the compactified moment map of the logarithmic cotangent bundle,
considered in \cite{Kn94}. Also, note that 
$R_{X_0} \cong G \times \fg$ over $X_0 \cong G$, and the restriction 
$f_0 : R_{X_0} \to \fg \times \fg$ is identified with the map
$$
G \times \fg \lto \fg \times \fg, \quad
(g,\xi) \longmapsto (\xi, g \cdot \xi).
$$ 
Thus, if $\xi \in \fg$ is regular and semi-simple, then the
fibre of $f_0$ at $(g,\xi)$ is isomorphic to the maximal torus 
$G^{g \cdot \xi}$. As a consequence, the general fibres of $f$ are
exactly the closures in $X$ of maximal tori of $G$.

Another natural map associates with any $x \in X$ the Lagrangian
subalgebra $\fg_{(x)}\subset \fg \times \fg$. In fact, this yields a
morphism from $X$ onto an irreducible component of the variety of
Lagrangian subalgebras, isomorphic to the wonderful compactification
of the adjoint semi-simple group $G/C(G)$ (see \cite[Sec.~2]{EL06}).
\end{example}

Returning to the general setting, we denote by $r(X)$ the dimension of
the general fibres of $f$; this is also the codimension of
$\overline{G \cdot \fh}$ in $\fg$, since $\dim R_X = \dim \fg$. 

\begin{corollary}\label{cor:ineq}
With the preceding notation and assumptions, we have
$r(X) \leq \rk(X)$.
\end{corollary}

\begin{proof}
We argue by induction on the rank of $X$. If $\rk(X) = 0$ then
$X \cong G/H$, where $H$ is a parabolic subgroup of $G$; thus, 
$\rk(H) = \rk(G)$. For an arbitrary rank $r$, consider a boundary
divisor $D_1$ and its open orbit $G \cdot x_1 \cong G/H_1$. Since
$\rk(D_1) = \rk(X) -1$, it suffices to show that 
\begin{equation}\label{eqn:rks}
\rk(H_1) \leq \rk(H) + 1.
\end{equation}
Choose a maximal torus $T_1 \subset H_1$ and consider its action on
the normal space to $D_1$ at $x_1$. This one-dimensional
representation defines a non-trivial character $\chi$ of $H_1$ (see
\cite[Prop.~2.1.2]{Bri07}) and hence of $T_1$. Let 
$S := \Ker(\chi\vert_{T_1})^0$. Then $S$ fixes points of $G/H$, as
follows from the existence of an \'etale linearisation of the action
of $T_1$ at $x_1$ (that is, a $T_1$-stable affine open subset
$X_1\subset X$ containing $x_1$, and  a $T_1$-equivariant \'etale
morphism $X_1 \to T_{x_1} X$ that maps $x_1$ to $0$). Thus, 
$\dim(S) \leq \rk(H)$, which yields (\ref{eqn:rks}). 
\end{proof}

\section{Vanishing theorems} 
\label{sec:vt}

\subsection{The linear case}
\label{subsec:lc}
We still consider a homogeneous pair $(X,D)$ under a connected linear
algebraic group $G$.

Since $X$ contains only finitely many orbits of $G$, Theorem
\ref{thm:finite} yields the vanishing of 
$H^i\big(X, \Omega_X^j(\log D)\big)$ for all $i > j$. From that
vanishing theorem for exterior powers of $\Omega_X^1(\log D)$, we
shall deduce a vanishing theorem for symmetric powers of
$\cR_X^{\vee}$, via the following homological trick (a generalisation
of \cite[Prop.~1]{Bro97}):

\begin{lemma}\label{lem:koszul}
Let $Z$ be a variety, and 
\begin{equation}\label{eqn:evf}
0 \lto \cE \lto \cO_Z \otimes_{\bC} V \lto \cF \lto 0
\end{equation}
an exact sequence of locally free sheaves, where $V$ is a 
finite-dimensional complex vector space. Then the following assertions
are equivalent for an invertible sheaf $\cL$ on $Z$ and an integer
$m$:

\smallskip

\noindent
{\rm (i)} $H^i(Z, \cL \otimes S^j \cF) = 0$ for all $i > m$ and all $j$.

\smallskip

\noindent
{\rm (ii)} $H^i(Z, \cL \otimes \wedge^j \cE) = 0$ for all $i > j + m$. 
\end{lemma}

\begin{proof}
Taking the Koszul complex associated with (\ref{eqn:evf}) and
tensoring with $\cL$ yields an exact sequence
$$\displaylines{
0 \to \cL \otimes \wedge^j \cE \to 
\cL \otimes \wedge^{j-1} \cE \otimes_{\bC} V \to \cdots 
\to \cL \otimes \wedge^{j-k} \cE \otimes_{\bC} S^k V \to \cdots \cr
\hfill \cdots \lto \cL \otimes \cE \otimes_{\bC} S^{j-1} V 
\to \cL \otimes S^j V 
\to \cL \otimes S^j \cF \lto 0.  
\cr}$$
We break this long exact sequence into short exact sequences:
$$
0 \lto \cL \otimes \wedge^j \cE \lto 
\cL \otimes \wedge^{j-1} \cE \otimes_{\bC} V \lto \cF_1 \lto 0,
$$
$$
0 \lto \cF_1 \lto \cL \otimes \wedge^{j-2} \cE \otimes_{\bC} S^2 V
\lto \cF_2 \lto 0, \quad \ldots
$$
$$
0 \lto \cF_{k-1} \lto \cL \otimes \wedge^{j-k} \cE \otimes_{\bC} S^k V
\lto \cF_k \lto 0, \quad \ldots
$$
$$
\quad 0 \lto \cF_{j-2} \lto \cL \otimes \cE \otimes_{\bC} S^{j-1}V
\lto \cF_{j-1} \lto 0,
$$
$$
0 \lto \cF_{j-1} \lto \cL \otimes S^j V \lto \cL \otimes S^j \cF \lto 0.
$$
If (i) holds, then $H^i(Z, \cL \otimes_{\bC} S^j V) = 0$ for all 
$i > m$, and hence $H^i(Z, \cF_{j-1}) = 0$ for all $i > m + 1$. We now
prove (ii) by induction on $j$. If $j = 1$, then $\cF_{j-1} = \cE$
which yields the assertion. For an arbitrary $j$, the induction
assumption implies that 
$H^i(Z, \cL \otimes \wedge^{j - k} \cE \otimes_{\bC} S^k V) = 0$ for all 
$k < j$ and $i > j - k + m$. By a decreasing induction on $k$, it
follows that $H^i(Z, \cF_k) = 0$ for all $i > j- k + m$. In particular, 
$ H^i(Z, \cF_1) = 0 
= H^i(Z, \cL \otimes \wedge^{j - 1} \cE \otimes_{\bC} V)$
for all $i >  j + m -1$, which implies the desired vanishing. 

The converse implication is obtained by reversing these arguments.
\end{proof}

We apply Lemma \ref{lem:koszul} to the exact sequence (\ref{eqn:dual})
and to $\cL = \cO_X$. Then the assertion (ii) holds for
$m = 0$; this yields:

\begin{theorem}\label{thm:sv}
With the assumptions of this subsection, we have 
\begin{equation}\label{eqn:van}
H^i(X, S^j \cR_X^{\vee}) = 0 \qquad (i \geq 1, ~j \geq 0).
\end{equation}
\end{theorem}

We now derive several geometric consequences of this vanishing
theorem. We shall need the following observation:

\begin{lemma}\label{lem:can}
The canonical sheaf of the nonsingular variety $R_X$ is given by
$\omega_{R_X} = p^*\cO_X(-D)$.
\end{lemma}

\begin{proof}
For any locally free sheaf $\cE$ of rank $r$ on $X$, the associated
vector bundle $E$ satisfies
$$
\omega_E = p^*(\wedge^r \cE^{\vee} \otimes \omega_X).
$$
Here, the rank of $\cR_X$ equals $\dim(G) - n$. Moreover,
$$
\wedge^{\dim(G) - n} \cR_X^{\vee} \cong \wedge^n \cT_X(-\log D) 
\cong \omega_X^{-1}(-D),
$$
as follows from (\ref{eqn:dual}) and (\ref{eqn:wed}).
\end{proof}

\begin{proposition}\label{prop:rat}
Denote by 
$$
\CD 
R_X @>{g}>> I_X @>{h}>> \fg
\endCD
$$
the Stein factorisation of the proper morphism $f$. Then $I_X$ is an
affine variety with rational singularities, and its canonical sheaf
satisfies
$$
\omega_{I_X} \cong R^{r(X)} f_* \big(p^*\cO_X(-D)\big).
$$
\end{proposition}

\begin{proof}
Since the morphism $h$ is finite, $I_X$ is affine; it is also normal,
since $R_X$ is nonsingular and the natural map 
$\cO_{I_X} \to g_* \cO_{R_X}$ is an isomorphism.

We claim that $R^i g_* \cO_{R_X} = 0$ for each $i \geq 1$. For this,
it suffices to show the vanishing of 
$h_*(R^i g_* \cO_{R_X}) = R^i f_* \cO_{R_X}$. Since the image of $f$
is affine, it suffices in turn to show that $H^i(R_X, \cO_{R_X}) = 0$. 
But this follows from (\ref{eqn:van}), since
$$
H^i(R_X, \cO_{R_X}) = H^i(X, p_* \cO_{R_X}) = 
\bigoplus_{j \geq 0} H^i(X, S^j \cR_X^{\vee}).
$$

We now deduce the rationality of singularities of $I_X$ from a result
of Koll\'ar (see \cite[Thm.~7.1]{Ko86}): let $\pi : Y \to Z$ be a
morphism of projective varieties, where $Y$ is nonsingular. If 
$\pi_* \cO_Y = \cO_Z$ and $R^i \pi_* \cO_Y = 0$ for all $i \geq 1$,
then $Z$ has rational singularities.

To reduce to that setting, we compactify the morphism $g$ as
follows. Consider the vector bundle $R_X \oplus O_X$ over $X$, where
$O_X$ denotes the trivial line bundle. 
This yields a subvariety of $X \times (\fg \oplus \bC)$, and hence a
proper morphism 
$$
\phi : R_X \oplus O_X \lto \fg \oplus \bC.
$$ 
By the preceding argument, 
\begin{equation}\label{eqn:vani}
R^i \phi_*\cO_{R_X \oplus O_X} = 0 \qquad (i \geq 1).
\end{equation}
Moreover, the set-theoretic fibre of $\phi$ at the origin of
$\fg \oplus \bC$ is the zero section. Thus, $\phi$ yields a morphism 
between projectivisations
$$
\bar{f}: \bP(R_X \oplus O_X) \lto \bP(\fg \oplus \bC)
$$
which extends $f: R_X \to \fg$. Furthermore, (\ref{eqn:vani}) easily
yields the vanishing of $R^i \bar{f}_*\cO_{\bP(R_X \oplus O_X)}$ for
$i \geq 1$. It follows that the Stein factorisation of $\bar{f}$,
$$
\bar{g}: \bP(R_X \oplus O_X) \lto \bar{I}_X
$$
extends $g$ and satisfies the same vanishing properties.

If the variety $X$ is projective, then so is $\bP(R_X \oplus O_X)$;
thus, Koll\'ar's result applies to $\bar{g}$ and hence to $g$. For an
arbitrary (complete, nonsingular) variety $X$, there exists a
nonsingular projective variety $X'$ together with a birational morphism
$$
u : X' \lto X.
$$ 
Then the pull-back to $X'$ of the projective bundle
$\bP(R_X \oplus O_X)$ is a projective variety $Y$ equipped with a
birational morphism 
$$
v : Y \lto \bP(R_X \oplus O_X).
$$ 
Since $\bP(R_X \oplus O_X)$ is nonsingular, we have 
$v_* \cO_Y = \cO_{\bP(R_X \oplus O_X)}$ and $R^i v_* \cO_Y = 0$ for 
all $i \geq 1$. Thus, Koll\'ar's result applies to the composite
morphism $\bar{g} \circ v$.

This completes the proof of the rationality of singularities of
$I_X$. The formula for its canonical sheaf follows from 
\cite[Theorem 5]{Ke76} in view of Lemma \ref{lem:can}.
\end{proof}

Next, we determine the algebra 
$$
H^{\bullet}\big(X, \Omega_X^{\bullet}(\log D)\big) :=
\bigoplus_{i,j} H^i\big(X, \Omega_X^j(\log D)\big)
$$
in terms of the coordinate ring of the affine variety $I_X$,  
\begin{equation}\label{eqn:coord}
\bC[I_X] = H^0(R_X,\cO_{R_X}) 
= \bigoplus_j  H^0(X, S^j \cR_X^{\vee}),
\end{equation}
viewed as a graded module over the algebra $\bC[\fg]$ (the symmetric
algebra of $\fg^*$) via the natural map 
$\fg^* \to H^0(X, \cR_X^{\vee})$. For this, consider the Koszul
complex associated with the exact sequence (\ref{eqn:dual}):
$$
\cdots \lto S^{\bullet} \cR_X^{\vee}\otimes_{\bC} \wedge^2 \fg^* 
\lto S^{\bullet} \cR_X^{\vee}\otimes_{\bC} \fg^*
\lto S^{\bullet} \cR_X^{\vee}.
$$
This complex of graded sheaves decomposes as a direct sum of complexes
$$
\cO_X \otimes_{\bC} \wedge^j \fg^* \lto 
\cR_X^{\vee} \otimes_{\bC} \wedge^{j-1} \fg^* \lto \cdots
\lto S^{j-1} \cR_X^{\vee} \otimes_{\bC} \fg^* \lto S^j \cR_X^{\vee}
$$
with homology sheaves $\Omega_X^j( \log D)$ in degree $0$, and $0$ in
positive degrees. Moreover, each sheaf 
$S^{j-k} \cR_X^{\vee} \otimes_{\bC} \wedge^k \fg^*$ is acyclic by
Theorem \ref{thm:sv}. This yields: 

\begin{proposition}\label{prop:comp}
With the assumptions of this subsection, each group
$H^i\big(\Omega_X^j( \log D)\big)$ is the $i$-th homology group of
the complex
$$\displaylines{
\wedge^j \fg^* \lto
H^0(X,\cR_X^{\vee}) \otimes \wedge^{j-1} \fg^* 
\lto \cdots 
\hfill \cr \hfill
\cdots \lto H^0(X, S^{j-1} \cR_X^{\vee}) \otimes \fg^* 
\lto H^0(X, S^j \cR_X^{\vee}).
\cr}$$
Moreover, we have isomorphisms
\begin{equation}\label{eqn:tor}
H^i\big(X, \Omega_X^j(\log D)\big) \cong 
\Tor_{j-i}^{j,\bC[\fg]}\big(\bC, \bC[I_X]\big),
\end{equation}
where $\bC$ is the quotient of $\bC[\fg]$ by its maximal
graded ideal, and the exponent $j$ denotes the subspace of degree $j$.
These isomorphisms are compatible with the multiplicative structures
of $H^{\bullet}\big(X, \Omega_X^{\bullet}(\log D)\big)$ and of 
$\Tor_{\bullet}^{\bullet,\bC[\fg]}\big(\bC, \bC[I_X]\big)$.
\end{proposition}

In turn, this will imply a simple description of the graded subalgebra
$H^0\big(X, \Omega_X^{\bullet}(\log D)\big)$.
To state it, consider the group $\cX(G)$ of multiplicative characters
of $G$, and its subgroup $\cX(G)^H$ of characters which restrict
trivially to $H$. Then $\cX(G)^H$ is a free abelian group of finite
rank, and every $f \in \cX(G)^H$ may be regarded as an invertible regular
function on $X_0 = G/H$; this yields an isomorphism 
$$
\cX(G)^H \cong \cO(X_0)^{\times}/\bC^{\times}.
$$
Also, $\cO(X_0)^{\times}/\bC^{\times}$ may be identified to
a subgroup of $H^0\big(X, \Omega_X^1(\log D)\big)$ via the map
$f \mapsto \frac{df}{f}$.

\begin{corollary}\label{cor:zero}
With the preceding notation, $H^0\big(X, \Omega_X^{\bullet}(\log D)\big)$
is a free exterior algebra on $\frac{df_1}{f_1}, \ldots,
\frac{df_r}{f_r}$, where $f_1,\ldots, f_r$ is any basis of the abelian
group $\cX(G)^H$.
\end{corollary}

\begin{proof}
Denote by $K$ the kernel of the map 
$\fg^* \to H^0(X, \cR_X^{\vee})$. Then $\wedge^j K$ is the kernel of
the induced map 
$\wedge^j \fg^* \to H^0(X, \cR_X^{\vee}) \otimes \wedge^{j-1} \fg^*$. 
By Proposition \ref{prop:comp}, it follows that 
$$
H^0\big(X, \Omega_X^{\bullet}(\log D)\big) \cong \wedge^{\bullet} K
$$
as graded algebras. In particular, 
$H^0\big(X, \Omega_X^1(\log D)\big) \cong K$. 

On the other hand, the exact sequence (\ref{eqn:res1}) yields an exact
sequence 
$$
0 \to H^0\big(X, \Omega_X^1(\log D)\big) 
\to \bC^{\ell} \lto H^1(X, \Omega_X^1) \to
H^1\big(\Omega_X^1(\log D)\big)  \to 0
$$
in view of Theorem \ref{thm:finite}. Moreover, the map 
$\bC^{\ell} \to H^1(\Omega_X^1)$ 
may be identified with the natural map
$\bigoplus_{k=1}^{\ell} A^0(D_k) \to A^1(X)$ 
as in the proof of that theorem. The kernel of the latter map is the
complexification of the abelian group 
$\{ \divi(f), ~f \in \cO(X_0)^{\times} \}$, and each $\divi(f)$ is
mapped to $\frac{df}{f}$ under the preceding identification.
\end{proof}

\begin{example}\label{exa:gr2}
Consider a $G \times G$-equivariant compactification $X$ of a
connected reductive group $G$, as in Example \ref{exa:gr1}. 
Then the image of $f : R_X \to \fg \times \fg$ is the closure of the 
set $\{(\xi, g \cdot \xi) ~\vert~ \xi \in \fg, ~g \in G\}$,
the graph of the adjoint action of $G$ on $\fg$. 

We claim that 
\begin{equation}\label{eqn:img}
I_X = f(R_X) = \fg \times_{\fg/\!/ G} \, \fg.
\end{equation}
Indeed, $f(R_X)$ is a variety of dimension $2 \dim(G) - \rk(G)$,
contained in $\fg \times_{\fg/\!/ G} \fg$. Moreover, since the
quotient morphism $\pi : \fg \to \fg/\!/G$ is flat and its
scheme-theoretic fibres are varieties of dimension $\dim(G) - \rk(G)$,
the same holds for the first projection 
$p_1 : \fg \times_{\fg/\!/ G} \, \fg \to \fg$. It follows that 
$\fg \times_{\fg/\!/ G} \, \fg$ is a variety of dimension  
$2 \dim (G) - \rk(G)$, which implies the second equality in 
(\ref{eqn:img}). To prove the first equality, it suffices to show that 
$\fg \times_{\fg/\!/ G} \, \fg$ is normal, since $R_X$ is smooth and the
general fibres of $f$ are connected. But $\fg \times_{\fg/\!/ G} \, \fg$
is a complete intersection in the affine space $\fg \times \fg$,
defined by the equations
$$
P_1(x) = P_1(y), \ldots, P_r(x) = P_r(y)
$$ 
where $P_1,\ldots,P_r$ are homogeneous generators of the graded
polynomial ring $\bC[\fg]^G$, and $r = \rk(G)$. In particular, 
$\fg \times_{\fg/\!/ G} \, \fg$ is Cohen-Macaulay. Moreover, the
differentials of $P_1,\ldots,P_r$ are linearly independent at any
regular element of $\fg$, and these form an open subset with
complement of codimension $3$. It follows that 
$\fg \times_{\fg/\!/ G} \, \fg$ is regular in codimension $1$, and hence
normal by Serre's criterion.

Together with (\ref{eqn:grp}) and (\ref{eqn:coord}), the equalities
(\ref{eqn:img}) imply that
$$
H^0\big(X, S^{\bullet} \cT_X(-\log D)\big) \cong 
S^{\bullet}(\fg) \otimes_{S^{\bullet}(\fg)^G} S^{\bullet}(\fg).
$$
Moreover, Theorem \ref{thm:sv} yields the vanishing of 
$H^i\big(X, S^{\bullet} \cT_X(-\log D)\big)$ for all $i \geq 1$, which
is also a special case of \cite[Thm.~4.1]{Kn94}. Combining both
results and arguing as in [loc.~cit.], we obtain an isomorphism
$$
\fU(X) \cong U(\fg) \otimes_{Z(\fg)} U(\fg),  
$$
where $\fU(X)$ denotes the algebra of differential operators on $X$
which preserve the ideal sheaf of $D$ (the ``completely regular''
differential operators of \cite{Kn94}), and $U(\fg)$ stands for the
enveloping algebra of $\fg$, with centre $Z(\fg)$.

On the other hand, (\ref{eqn:tor}) yields isomorphisms
$$
H^i\big(X,\Omega_X^j(\log D)\big) \cong 
\Tor_{j-i}^{j,\bC[\fg] \otimes \bC[\fg]}
\big(\bC, \bC[\fg] \otimes_{\bC[\fg]^G} \bC[\fg]\big).
$$
Since the algebra $\bC[\fg] \otimes_{\bC[\fg]^G}\bC[\fg]$ is the
quotient of the polynomial algebra $\bC[\fg \times \fg]$ by the
ideal generated by the regular sequence
$\big(P_1(x)- P_1(y), \ldots,P_r(x) - P_r(y)\big)$, 
this yields in turn
$$
H^i\big(X,\Omega_X^j(\log D)\big) \cong 
\Tor_{j-i}^{j,\bC[\fg]^G}(\bC, \bC).
$$
As a consequence, the bi-graded algebra 
$H^{\bullet}\big(X,\Omega_X^{\bullet}(\log D)\big)$
is a free exterior algebra on generators of bi-degrees
$(d_1 - 1,d_1), \ldots, (d_r - 1,d_r)$, where $d_1,\ldots,d_r$ denote
the degrees of $P_1,\ldots,P_r$. 
\end{example}

\begin{remark}
More generally, consider a log homogeneous variety $X$ under a
connected reductive group $G$, and assume that $H$ is connected and
reductive as well. Then the invariant rings $\bC[\fg]^G$ and
$\bC[\fh]^H$ are graded polynomial rings, and $\bC[\fh]^H$ is a finite
module over $\bC[\fg]^G$ via the restriction map. Moreover, the
quotient morphism $\fh \to \fh/\!/H  := \Spec \bC[\fh]^H$ yields a
morphism $G \times^H \fh \to \fh/\!/H$, the quotient of the affine
variety $G \times^H \fh$ by the action of $G$. 
One may check that the latter morphism extends to $R_X$, and that the 
product map $R_X \to  \fg \times \fh/\!/H$ factors through an
isomorphism
$$
I_X \cong \fg \times_{\fg/\!/G} \fh/\!/H.
$$
In other words,
$$
\bC[I_X] \cong \bC[\fg] \otimes_{\bC[\fg]^G} \bC[\fh]^H.
$$
Together with the isomorphism (\ref{eqn:tor}) and the freeness of the
$\bC[\fg]^G$-module $\bC[\fg]$, it follows that
$$
H^i\big( X, \Omega_X^j(\log D) \big) \cong 
\Tor_{j-i}^{j,\bC[\fg]^G}\big( \bC, \bC[\fh]^H \big).
$$
Also, note that $\bC[\fg]^G$ (resp. $\bC[\fh]^H$) is the cohomology
ring of the classifying space $BG$ (resp. $BH$) with complex
coefficients. This yields a description of the algebra 
$H^{\bullet}\big( X, \Omega_X^{\bullet}(\log D) \big)$ in topological
terms, which can be extended to any homogeneous space -- not
necessarily having a log homogeneous compactification -- in view of a
result of Franz and Weber (see \cite[Thm.~1.6]{FW05}).
\end{remark}

\subsection{The linear case (continued)}
\label{subsec:lcc}

We still consider a homogeneous $G$-pair $(X,D)$, where $G$ is a
connected linear algebraic group. Let $\cL$ denote an invertible sheaf
on $X$, and assume that $\cL$ is nef; since $X$ is a spherical
variety under a Levi subgroup of $G$, this is equivalent to $\cL$ being  
generated by its global sections (see e.g. \cite[Lem.~3.1]{Bri93}).  

Recall from Subsec.~\ref{subsec:bils} that any component $F$ of a
general fibre of $f: R_X \to \overline{G \cdot \fh}$ may be
identified to a toric subvariety of $X$. 
We denote by $\kappa_f(\cL)$ the Kodaira--Iitaka dimension of the
pull-back $\cL \vert_F$. Since $\cL$ is globally generated,
$\kappa_f(\cL)$ is the dimension of the image of the natural map 
$\varphi : F \to \bP\big(H^0(F,\cL)^*\big)$.
Note that $\kappa_f(\cL) \leq \dim(F) = r(X)$, and equality holds
e.g. if $\cL$ is ample.

\begin{theorem}\label{thm:van}
With the notation and assumptions of this subsection,
$H^i\big(X, \cL(-D) \otimes S^j \cR_X^{\vee} \big) = 0$
for all $i \neq r(X) - \kappa_f(\cL)$ and $j \geq 0$. 
\end{theorem}

\begin{proof}
Since $R_X$ is nonsingular, $f$ is proper and the invertible sheaf
$p^*\cL$ is $f$-semi-ample, it follows that the sheaf 
$R^i f_*(p^* \cL \otimes \omega_{R_X})$ is torsion-free on the image
of $f$, for any $i \geq 0$ (see \cite[Cor.~6.12]{EV92}). But
$$
H^i(F, \cL \otimes \omega_{R_X}) =  H^i(F, \cL \otimes \omega_F)
$$
vanishes for all $i \neq r(X) - \kappa_f(\cL)$, by a result of Fujino
(see \cite[Cor.~1.7]{Fuj07}). Thus, the same vanishing holds for 
$R^i f_*(p^* \cL \otimes \omega_{R_X})$, 
i.e., for $R^i f_*\big(p^* \cL(-D)\big)$ in view of Lemma
\ref{lem:can}. This implies our statement, by arguing as in 
the beginning of the proof of Proposition~\ref{prop:rat}.
\end{proof}

In view of Lemma \ref{lem:koszul}, Theorem \ref{thm:van} yields
the vanishing of the groups
$H^i\big(X, \cL(-D) \otimes \Omega_X^j(\log D)\big)$ 
for all $i > j + r(X) - \kappa_f(\cL)$. By Serre duality
(\ref{eqn:sd}), this implies the following: 

\begin{corollary}\label{cor:van}
With the notation and assumptions of this subsection, 
$H^i\big(X, \cL^{-1} \otimes \Omega_X^j(\log D)\big) = 0$ 
for all $i < j - r(X) + \kappa_f(\cL)$.

In particular, this vanishing holds for all $i < j$ if $\cL$ is
ample. 
\end{corollary}

\begin{remarks}
(i) Assume that $H^k(F,\cL^{-1}) \neq 0$ for some $k\geq 0$ and some
nef invertible sheaf $\cL$ on $X$. Then $k = \kappa_f(\cL)$, and
$$
H^{r(X) -k}\big(X, \cL(-D) \otimes S^j \cR_X^{\vee} \big) \neq 0
$$
for some $j \geq 0$, by the proof of Theorem \ref{thm:van}. 
Equivalently, 
$$
H^i\big(X, \cL^{-1} \otimes \Omega_X^j(\log D)\big) \neq 0
$$ 
for some $i = j - r(X) + k$. 

The preceding assumption on $\cL$ is fulfilled if $\cL = \cO_X$; then
$k = 0$. Thus, Theorem \ref{thm:van} and Corollary \ref{cor:van} 
are optimal for the trivial invertible sheaf. 

The non-vanishing of $H^{r(X)}\big(X, (S^j \cR_X^{\vee})(-D)\big)$
for some $j$ also follows from Proposition \ref{prop:rat}: since $I_X$
is affine, the space of global sections of $\omega_{I_X}$ is non-zero.

Likewise, Theorem \ref{thm:van} and Corollary \ref{cor:van} are
optimal for ample invertible sheaves, since 
$H^0\big(F,\cL(-D)\big) \neq 0$ for sufficiently ample $\cL$;
equivalently, $H^{r(X)}(F,\cL^{-1}) \neq 0$.

\smallskip

\noindent
(ii) If $\cL$ is big, then we also have a refinement of Norimatsu's
vanishing theorem mentioned in Subsec.~\ref{subsec:dflp}; namely,
$$
H^i\big(X, \cL^{-1} \otimes \Omega_X^j(\log D)\big) = 0 \qquad 
(i + j < n),
$$ 
as follows from \cite[Cor.~6.7]{EV92}. 

\end{remarks}

\subsection{The case where the open orbit is proper over an affine}
\label{subsec:aoe}

We consider again a homogeneous $G$-pair $(X,D)$, where $G$ is a
connected linear algebraic group, and a nef invertible sheaf $\cL$ on
$X$. We assume in addition that the open orbit $X_0$ is proper over an
affine; this holds e.g. if $H$ is reductive (and hence $X_0$ is
affine), or if $X$ is a flag variety. We now obtain a stronger version
of the vanishing theorem \ref{thm:sv}: 

\begin{theorem}\label{thm:aff}
With the assumptions of this subsection, we have 
$H^i(X, \cL \otimes S^j \cR_X^{\vee}) = 0$ 
for all $i \geq 1$ and all $j$. 

Equivalently, $H^i\big(X, \cL \otimes \Omega_X^j(\log D)\big) = 0$ 
for all $i > j$.
\end{theorem}

\begin{proof}
By Lemma \ref{lem:ind}, we may replace $X$ with $X'$ and
$\cL$ with $\cL' := u^* \cL$, where $X'$ is a complete nonsingular
$G$-variety and $u: X' \to X$ is a $G$-equivariant birational
morphism. We claim that we may choose $X'$ so that $D'$
(the reduced inverse image of $D$) is the support of an effective
base-point-free divisor.

By assumption, we have a proper morphism $\varphi : X_0 \to Y_0$,
where $Y_0$ is an affine variety. We may assume in addition that 
$\varphi_*\cO_{X_0} = \cO_{Y_0}$. Then $\varphi$ is a surjective
$G$-equivariant morphism with connected fibres, and hence a fibration
in flag varieties. The affine $G$-variety $Y_0$ admits a closed
$G$-equivariant immersion into some $G$-module $V$. Consider the 
associated rational map $X - \to \bP(V \oplus \bC)$, and its graph
$X'$. Then $X'$ is a complete $G$-variety, and the projection 
$u : X' \to X$ is an isomorphism above $X_0$. Moreover, the complement
$X' \setminus u^{-1}(X_0) = \Supp\big( u^{-1}(D)\big)$ is the
set-theoretic preimage of the hyperplane section $\bP(V \oplus 0)$
under the projection $X' \to \bP(V\oplus \bC)$. Now replace $X'$ with
an equivariant desingularisation to obtain the desired setting.

Thus, we may assume that there exists a base-point-free divisor
$\sum_{k=1}^{\ell} a_k \, D_k$, where the $a_k$ are positive 
integers. Choose an integer $N > a_1, \ldots, a_{\ell}$. We now apply
\cite[Cor.~6.12]{EV92} to the morphism 
$f : R_X \to \overline{G \cdot \fh}$, the invertible sheaf 
$\cM := p^* \cL(D)$, and the divisor 
$E := p^*\big(\sum_{k=1}^{\ell} (N - a_k) \, D_k\big)$.
Then $E$ has normal crossings, and the invertible sheaf
$$
\cM^N(-E) = p^* \cL^N(\sum_{k=1}^{\ell} a_k \, D_k)
$$
is $f$-semi-ample, so that the assumptions of [loc.~cit.] are
satisfied. Hence each sheaf $R^i f_*(p^* \cM \otimes \omega_{R_X})$ is
torsion-free on the image of $f$. By Lemma \ref{lem:can}, this means
that $R^i f_*(p^* \cL)$ is torsion-free. But $H^i(F, \cL)=0$ for
any $i \geq 1$ and any component $F$ of a general fibre of $f$, in
view of Theorem \ref{thm:fibres}. Thus, $R^i f_*(p^* \cL) = 0$ for all 
$i \geq 1$. This implies our statements by the arguments of the
preceding subsection.
\end{proof}

\begin{remark}
The preceding argument also yields a simpler proof of
\cite[Thm.~3.2]{BB96}, the main result of that paper. It asserts that 
$H^i\big(X, \cL \otimes S^j\cT_X(- \log D)\big) = 0$ for all 
$i \geq 1$ and $j\geq 0$, where $(X,D)$ is a homogeneous pair under a
connected reductive group $G$, the open orbit $X_0$ is proper over an
affine, and $\cL$ is a nef invertible sheaf on $X$. Here the morphism 
$f : R_X \to \fg$ is replaced with the compactified moment map of
\cite{Kn94}. 
\end{remark}

\subsection{The general case}
\label{subsec:gc}

We now consider a $G$-pair $(X,D)$, where $G$ is a connected algebraic
group, not necessarily linear. 

We shall obtain a generalisation of Corollary \ref{cor:van} to this
setting. Note that the arguments of Subsec.~\ref{subsec:lcc} need to
be substantially modified, since the assumption of semi-ampleness in
Koll\'ar's result is not satisfied (e.g., for algebraically trivial
invertible sheaves on abelian varieties). Thus, we begin with a closer
study of the Albanese fibration.

With the notation of Subsec.~\ref{subsec:vfmog}, we identify $X$ to 
$G \times^I Y$, and the Albanese morphism  $\alpha : X \to A$ to the
natural map $G \times^I Y \to G/I$. Let $E := D \vert_Y$, then $(Y,E)$
is a homogeneous pair under $G_{\aff} = I^0$, and $I$ preserves each
$G_{\aff}$-orbit in $Y$ (see \cite[Thm.~3.2.1]{Bri07}). As a
consequence, $\rk(X) = \rk(Y)$. Also, recall that $Y$ is a spherical
variety under any Levi subgroup of $G_{\aff}$.

We denote by $C(G)$ the centre of $G$, so that
$$
G = C(G)^0 G_{\aff}
$$ 
(see e.g. \cite[Lem.~1.1.1]{Bri07}). Thus, $C(G)^0$ acts transitively
on $A$, and $I = \big(I \cap C(G)^0\big)G_{\aff}$. Also, $C(G)^0$ is a
semi-abelian variety by \cite[Prop.~3.4.2]{Bri07}. Moreover, the
isotropy subgroup $G_x$ is affine for any $x \in X$ (see
e.g. \cite[Lem.~1.2.1]{Bri07}); as a consequence, 
$\fg_{(x)} \subset \fg_x \subset \fg_{\aff}$. It follows 
that
$$
R_X \cong G \times^I R_Y \quad \text{and} \quad 
\overline{G \cdot \fh} = \overline{G_{\aff} \cdot \fh} \subset
\fg_{\aff}. 
$$
Moreover, denoting by $X^{(\xi)}$ (resp.~$Y^{(\xi)}$) the fibre at 
$\xi \in \fg_{\aff}$ of the map $f$ (resp.~$f_Y: R_Y \to \fg_{\aff}$),
we see that $Y^{(\xi)}$ is preserved by $I$, and 
$$
X^{(\xi)} = G \times^I Y^{(\xi)}.
$$ 
In particular, \emph{the connected components of the general fibres of
$f$ are semi-abelic varieties of dimension $q(X) + r(Y)$.} In view of
Corollary \ref{cor:ineq}, this yields
\begin{equation}\label{eqn:rqr}
r(X) = q(X) + r(Y) \leq q(X) + \rk(X).
\end{equation}

Next, we describe the invertible sheaves on $X$. For this, choose a
Borel subgroup $B \subset G_{\aff}$. Then 
$$
G_1 := C(G)^0 B
$$ 
is a maximal connected solvable subgroup of $G$, and $(G_1)_{\aff} = B$.

\begin{lemma}\label{lem:prod}
{\rm (i)} Any invertible sheaf $\cL$ on $X$ admits a decomposition
\begin{equation}\label{eqn:dec}
\cL = (\alpha^*\cM)(\Delta), 
\end{equation}
where $\cM$ is an invertible sheaf on $A$, and $\Delta$ is a
$G_1$-stable divisor on $X$. Moreover, $\cM$ (resp.~$\Delta\vert_Y$)
is uniquely determined by $\cL$ up to algebraic (resp.~rational)
equivalence.

\smallskip

\noindent
{\rm (ii)} $\cL$ is nef (resp.~ample) if and only if both $\cM$ and
$\Delta \vert_Y$ are nef (resp.~ample).
\end{lemma}

\begin{proof}
(i) Note that $B$ has an open orbit $Y_1$ in $Y$, and hence 
$G_1$ has an open orbit $X_1$ in $X$; the map
$$
\alpha_1 := \alpha \vert_{X_1} : X_1 \lto A
$$
is a $G_1$-equivariant fibration with fiber $Y_1$.

We claim that the pull-back map $\alpha_1^* : \Pic(A) \to \Pic(X_1)$
is surjective. To see this, identify $X_1$ to the homogeneous space
$G_1/H_1$; then $A = X/B = G_1/H_1 B$ and hence $\alpha_1$ is the
composite morphism
$$
\CD
G_1/H_1 @>{\alpha_U}>> G_1/H_1 U @>{\alpha_T}>> G_1/H_1 B,
\endCD
$$ 
where $U$ denotes the unipotent part of $B$, and $T$ denotes the torus
$B/U$. Note that $G_1/U$ is a semi-abelian variety with maximal torus
$T$. Thus, the quotient $G_1/H_1 U$ is a semi-abelian variety as well,
and $\alpha_T$ is the quotient map by its maximal torus. It follows
that 
$$
\alpha_T^* : \Pic(A) \lto \Pic(G_1/H_1 U)
$$ 
is surjective. On the other hand, since $U$ is unipotent, $\alpha_U$
may be factored into quotients by free actions of the additive group,
and hence
$$
\alpha_U^* : \Pic(G_1/H_1 U) \lto \Pic(G_1/H_1)
$$ 
is an isomorphism.
 
By the claim, there exists an invertible sheaf $\cM$ on $A$
such that $\cL \vert_{X_1} \cong \alpha_1^* \cM$. Then 
$\cL \cong (\alpha^* \cM)(\Delta)$ for some divisor
$\Delta$ supported in $X \setminus X_1$. In particular, $\Delta$ is
preserved by $G_1$. This proves the existence of the decomposition
(\ref{eqn:dec}). 

For the uniqueness properties, we may assume that $\cL$ is
trivial. Then $\alpha^*(\cM) = \cO_X(- \Delta)$, and hence
$g^* \alpha^*(\cM) \cong \alpha^*(\cM)$ 
for any $g \in G_1$. It follows that $a^*(\cM) \cong \cM$
for any $a \in A$; thus, $\cM$ is algebraically trivial. Moreover,
$\cO_Y(\Delta\vert_Y) = \cL \otimes \cO_Y$ is trivial as well.

(ii) Recall that any effective $1$-cycle on $X$ is rationally
equivalent to an effective $1$-cycle preserved by $B$ (see
e.g. \cite[Sec.~1.3]{Bri93}). Thus, $\cL$ is nef if and only if 
$\cL \cdot C \geq 0$ for any irreducible curve $C \subset X$,
preserved by $B$. 

If $B$ acts non-trivially on $C$, then $C$ is rational and
hence contained in a fibre of $\alpha$. Thus, 
$\cL \cdot C = \Delta \cdot C$. Since $\Delta$ is preserved by the
group $B_1$ which permutes transitively the fibres of $\alpha$, we may
assume that $C \subset Y$; then 
$\cL \cdot C = \Delta\vert_Y \cdot C$.

On the other hand, if $B$ acts trivially on $C$, then the orbit 
$G \cdot x$ is closed in $X$ for any $x \in C$. Since $X$ contains
only finitely many $G$-orbits, it follows that $C \subset G \cdot x$
for any such $x$. By \cite[Thm.~3.3.3]{Bri07}, there is a
$G$-equivariant isomorphism 
$
G \cdot x \cong A \times (G_{\aff} \cdot x).
$
This identifies the fixed point subscheme 
$(G\cdot x)^B$ to $A \times \{x\}$, and hence $C$ to a curve in 
$A \times \{x\}$; in particular, 
$\alpha$ restricts to an isomorphism $C \cong \alpha(C)$. Since 
$A \times \{x\}$ is an orbit of $G_1$,
the restriction $\Delta \vert_{A \times \{x\}}$ is algebraically trivial
by the argument of (ii). It follows that
$\cL \cdot C = (\alpha^* \cM) \cdot C = \cM \cdot \alpha(C)$.

Thus, $\cL$ is nef iff so are $\cM$ and $\Delta \vert_Y$. This implies
the corresponding statement for ampleness, in view of Kleiman's
criterion: $\cL$ is ample iff for any invertible sheaf $\cL'$ on $X$,
there exists a positive integer $n = n(\cL')$ such that 
$\cL^n \otimes \cL'$ is nef. 
\end{proof}

\begin{remarks}\label{rem:prod}
(i) Lemma \ref{lem:prod} implies readily a decomposition of the 
N\'eron--Severi group:
$$
\NS(X) \cong \NS(A) \times \NS(Y).
$$
Also, $\NS(A)$ is a free abelian group; moreover, 
$\NS(Y)$ is also free and isomorphic to $\Pic(Y)$, since $Y$
admits a cellular decomposition (Lemma \ref{lem:pav}). It follows that
$\NS(X)$ is a free abelian group as well; in other words, algebraic
and numerical equivalence coincide for invertible sheaves on $X$. 
Moreover, the cone of numerical equivalence classes of nef invertible
sheaves decomposes accordingly:
$$
\Nef(X) \cong \Nef(A) \times \Nef(Y).
$$
The nef cones of abelian varieties are well understood (see
e.g. \cite{BL04}). Those of spherical varieties (like $Y$) are studied
in \cite{Bri93}; in particular, these cones are polyhedral. 

\smallskip

\noindent
{\rm (ii)} If $X$ is a semi-abelic variety, then $G_{\aff}$ is a torus
$T$, and $Y$ is a toric variety under that torus. Thus, $G_1 = G$, and
$G_1$-stable divisors on $X$ correspond bijectively to $T$-stable
divisors on $Y$. In that case, Lemma \ref{lem:prod} gives back a
description of the ample invertible sheaves on semi-abelic varieties,
due to Alexeev (see \cite[Sec.~5.2]{Al02}).
\end{remarks}

Consider a nef invertible sheaf $\cL$ on $X$, and its decomposition
(\ref{eqn:dec}); then $\cL\vert_Y \cong \cO_Y(\Delta \vert_Y)$. Let 
$$
K(\cM) := \{a \in A ~\vert~ a^*\cM \cong \cM\}; 
$$
this is a closed subgroup of $A$. Since $\cL$ determines $\cM$ up to
multiplication by an algebraically trivial invertible sheaf, we see
that $K(\cM)$ depends only on $\cL$; we shall denote that group by
$K(\cL)$. We now are in a position to state:

\begin{theorem}\label{thm:vang}
With the preceding notation and assumptions, we have 
$H^i\big(X,\cL^{-1} \otimes \Omega_X^j(\log D)\big)  = 0$ for all 
$i < j - r(Y) + \kappa_f(\cL\vert_Y) - \dim K(\cL)$. 

In particular, this vanishing holds for all $i < j$ if $\cL$ is ample.
\end{theorem}

\begin{proof}
By (\ref{eqn:sd}), it suffices to show that
$H^i\big(X,\cL(-D) \otimes \Omega_X^j(\log D)\big) = 0$
for all $i > j + r(Y) - \kappa_f(\cL\vert_Y) + \dim K(\cL)$.
In view of the decomposition (\ref{eqn:albj}), this reduces to showing
that 
\begin{equation}\label{eqn:gen}
H^i\big(X,\cL(-D) \otimes \Omega_{X/A}^j(\log D)\big) = 0
\end{equation}
for all such $i$ and $j$. Consider the Leray spectral sequence
$$\displaylines{
E_2^{p,q} := 
H^p\big(A, R^q \alpha_*(\cL(- D) \otimes \Omega_{X/A}^j(\log D))\big)
\hfill \cr \hfill
\Rightarrow 
H^{p+q}\big(X,\cL(- D) \otimes \Omega_{X/A}^j(\log D)\big).
\cr}$$
By the projection formula, 
$$
R^q \alpha_*\big(\cL(- D) \otimes \Omega_{X/A}^j(\log D)\big)
\cong \cM \otimes 
R^q \alpha_*\big(\Omega_{X/A}^j(\log D)(- \Delta - D)\big).
$$
Note that the sheaf 
$R^q \alpha_*\big(\Omega_{X/A}^j(\log D)(- \Delta - D)\big)$
is locally free and $G_1$-linearised. Hence this sheaf has a
filtration with subquotients being algebraically trivial invertible
sheaves. Since $\cM$ is nef, it follows that 
$E_2^{p,q} = 0$ for any $p > \dim K(\cL) = \dim K(\cM)$, by a
classical vanishing theorem for abelian varieties (see 
\cite[Lem.~3.3.1, Thm.~3.4.5]{BL04}). On the other hand,
since $\Delta\vert_Y$ is nef, Corollary \ref{cor:van} yields that 
$H^q\big(Y,\Omega_Y^j(\log E)(- \Delta\vert_Y - E)\big) = 0$
for any $q > i + r(Y) - \kappa_f(\cL\vert_Y)$, and hence 
$E_2^{p,q} = 0$ for all such $q$. This implies (\ref{eqn:gen}).
\end{proof}

Taking for $\cL$ the trivial invertible sheaf and combining 
Theorems \ref{thm:finiteg} and \ref{thm:vang} with the inequality
(\ref{eqn:rqr}), we obtain: 

\begin{corollary}\label{cor:strip}
With the notation and assumptions of this subsection, 
$H^i\big(X, \Omega_X^j(\log D)\big) = 0$ unless
$0 \leq j - i \leq r(X)$.
\end{corollary}

Another consequence of Theorem \ref{thm:vang} is a vanishing result
for ordinary Dolbeault cohomology:

\begin{theorem}\label{thm:dol}
Let $\cL$ be an invertible sheaf on a log homogeneous variety
$X$ of irregularity $q$ and rank $r$.

If $\cL$ is nef (resp.~ample), then 
$H^i\big(X, \cL^{-1} \otimes \Omega_X^j(\log D')\big) = 0$
for any effective subdivisor $D'$ of $D$, and for all
$i < j - q - r$ (resp.~$ i < j$). 

In particular, 
$H^i\big(X, \cL \otimes \Omega_X^j\big) = 0$ for all  
$i > j + q + r$ if $\cL$ is nef, resp.~ $i > j$ if $\cL$ is ample.
\end{theorem}

\begin{proof}
If $D' = D$, then the first assertion is a consequence of Theorem 
\ref{thm:vang} and Lemma \ref{lem:prod}. Indeed, if $\cL$ is nef, then
$r(Y) -\kappa_f(\cL\vert_Y) \leq \rk(Y) = r$ by Corollary
\ref{cor:ineq}, and $\dim K(\cL) \leq q$. If $\cL$ is ample, then we
have the equalities 
$r(Y) - \kappa_f(\cL\vert_Y) = 0 = \dim K(\cL)$.

The case of an arbitrary divisor $D'$ follows by decreasing induction
on the number of irreducible components of $D'$ and the dimension of
$X$, using (\ref{eqn:resj}).

Taking $D' = 0$ yields the second assertion, by Serre duality.
\end{proof}


\begin{thebibliography}{100} 


\bibitem[Al02]{Al02}
V.~Alexeev,
\emph{Complete moduli in the presence of semiabelian group actions},
Ann. of Math. (2) \textbf{155} (2002), no.~3, 611--708.


\bibitem[BC94]{BC94}
V.~V.~Batyrev and D.~A.~Cox,
\emph{On the Hodge structure of projective hypersurfaces in toric
varieties}, Duke Math. J. \textbf{75} (1994), no.~2, 293--338.


\bibitem[Bi73]{Bi73}
A.~Bia\l ynicki-Birula,
\emph{Some theorems on actions of algebraic groups}, 
Ann. of Math. (2) \textbf{98} (1973), 480--497. 


\bibitem[BB96]{BB96}
F.~Bien and M.~Brion,
\emph{Automorphisms and local rigidity of regular varieties},
Compositio Math. \textbf{104} (1996), 1--26.


\bibitem[BL04]{BL04} C.~Birkenhake and H.~Lange,
\emph{Complex abelian varieties. Second edition},
Springer--Verlag, Berlin, 2004.




\bibitem[Bri93]{Bri93}
M.~Brion,
\emph{Vari\'et\'es sph\'eriques et th\'eorie de Mori},
Duke Math. J. \textbf{72} (1993), no.~2, 369--404.


\bibitem[Bri07]{Bri07} 
M.~Brion,
\emph{Log homogeneous varieties},
Actas del XVI Coloquio Latinoamericano de \'Algebra, 1--39,
Revista Matem\'atica Iberoamericana, Madrid, 2007;
arXiv: math/0609669.


\bibitem[Bro93]{Bro93}
A.~Broer,
\emph{Line bundles on the cotangent bundle of the flag variety},
Invent. math. \textbf{113} (1993), no.~1, 1--20.


\bibitem[Bro97]{Bro97}
A.~Broer,
\emph{A vanishing theorem for Dolbeault cohomology of homogeneous 
vector bundles},
J. reine angew. Math. \textbf{493} (1997), 153--169.


\bibitem[CL73]{CL73}
J.~B.~Carrell and D.~I.~Lieberman,
\emph{Holomorphic vector fields and compact Kaehler manifolds},
Invent. math. \textbf{21} (1973), 303--309.


\bibitem[DK86]{DK86}
V.~I.~Danilov and A.~.G.~Khovanskii, 
\emph{Newton polyhedra and an algorithm for calculating Hodge--Deligne
numbers}, 
Math. USSR-Izv. \textbf{29} (1987), no.~2, 279--298.


\bibitem[DP83]{DP83}
C.~de Concini and C.~Procesi, 
\emph{Complete Symmetric Varieties}, 
Lect. Notes in Math., Springer, \textbf{996} (1983), 1--44.


\bibitem[De71]{De71}
P.~Deligne, 
\emph{Th\'eorie de Hodge II},
Pub. Math. IHES \textbf{40} (1971), 5--57.


\bibitem[De74]{De74}
P.~Deligne, 
\emph{Th\'eorie de Hodge III},
Pub. Math. IHES \textbf{44} (1974), 5--78.


\bibitem[EL06]{EL06}
S.~Evens and J.~H.~Liu,
\emph{On the variety of Lagrangian subalgebras II},
Ann. scient. \'Ecole Norm. Sup. (4) \textbf{39} (2006), no.~2,
347--379. 


\bibitem[EV92]{EV92} 
H.~Esnault and E.~Viehweg,
\emph{Lectures on Vanishing Theorems},
DMV Seminar \textbf{20}, Birkh\"auser, Basel, 1992.


\bibitem[FW05]{FW05}
M.~Franz and A.~Weber,
\emph{Weights in cohomology and the Eilenberg--Moore spectral sequence},
Ann. Inst. Fourier, Grenoble \textbf{55}, no.~2 (2005), 673--691.


\bibitem[Fuj07]{Fuj07}
O.~Fujino,
\emph{Multiplication map and vanishing theorems for toric varieties},
Math. Z. \textbf{257} (2007), no.~3, 631--641.


\bibitem[Ful98]{Ful98} W.~Fulton,
\emph{Intersection Theory. Second edition},
Ergeb. Math. Grenzgebiete. 3. Folge. Springer-Verlag, Berlin, 1998.


\bibitem[Ha77]{Ha77}
R.~Hartshorne,
\emph{Algebraic Geometry},
Graduate Texts in Mathematics \textbf{52}, Springer-Verlag, New York,
1977. 


\bibitem[Ke76]{Ke76}
G.~Kempf, 
\emph{On the collapsing of homogeneous bundles},
Invent. math. \textbf{37} (1976), no.~3, 229--239.


\bibitem[Ki06]{Ki06}
V.~Kiritchenko,
\emph{Chern classes of reductive groups and an adjunction formula},
Ann. Inst. Fourier (Grenoble) \textbf{56} (2006), no.~4, 1225--1256.


\bibitem[Kn94]{Kn94}
F.~Knop,
\emph{A Harish-Chandra homomorphism for reductive group actions},
Ann. of Math. (2) \textbf{140} (1994), no.~2, 253--288.


\bibitem[Ko86]{Ko86}
J.~Koll\'ar,
\emph{Higher direct images of dualizing sheaves},
Ann. of Math. (2) \textbf{123} (1986), no.~1, 11-42.


\bibitem[Lu01]{Lu01}
D.~Luna,
\emph{Vari\'et\'es sph\'eriques de type $A$},
Pub. Math. IH\'ES \textbf{94} (2001), 161--226.


\bibitem[Ma99]{Ma99}
A.~V.~Mavlyutov, 
\emph{Cohomology of complete intersections in toric varieties},
Pacific J. Math. \textbf{191} (1999), no. 1, 133--144.


\bibitem[Ma08]{Ma08}
A.~V.~Mavlyutov, 
\emph{Cohomology of rational forms and a vanishing theorem for toric
varieties}, 
J. reine angew. Math. \textbf{615} (2008), 45--58. 


\bibitem[No78]{No78}
Y.~Norimatsu, 
\emph{Kodaira vanishing theorem and Chern classes for 
$\partial$-mani\-folds},
Proc. Japan Acad. Ser. A Math. Sci. \textbf{54} (1978), no.~4, 107--108.


\bibitem[Sn86]{Sn86}
D.~Snow,
\emph{Cohomology of twisted holomorphic forms on Grassmann manifolds
and quadric hypersurfaces}, Math. Ann. \textbf{276} (1986), no.~1,
159--176. 


\bibitem[To08]{To08}
B.~Totaro,
\emph{Chow groups, Chow cohomology, and linear varieties},
preprint available at www.dpmms.cam.ac.uk/\~{}bt219/papers.html, 
to appear in Journal of Algebraic Geometry.



\bibitem[We03]{We03}
J.~Weyman,
\emph{Cohomology of vector bundles and syzygies},
Cambridge Tracts in Mathematics \textbf{149}, Cambridge University Press,
Cambridge, 2003. 


\bibitem[Wi04]{Wi04}
J.~Winkelmann, 
\emph{On manifolds with trivial logarithmic tangent bundle},
Osaka J. Math. \textbf{41} (2004), no. 2, 473--484. 

\end{thebibliography}
\end{document}